\documentclass[11pt]{extarticle}
\usepackage[english]{babel}
\usepackage{graphicx}
\usepackage{harvard}
\usepackage{framed}
\usepackage[normalem]{ulem}
\usepackage{amsmath}
\usepackage{amsthm}
\usepackage{amssymb}
\usepackage{amsfonts}
\usepackage{enumerate}
\usepackage{xcolor}
\usepackage[utf8]{inputenc}
\usepackage{dsfont}
\usepackage{comment}
\usepackage[top=1 in,bottom=1in, left=1 in, right=1 in]{geometry}
\usepackage{hyperref}

\newcommand{\cov}{{\rm Cov}}
\newcommand{\var}{{\rm Var}}

\newcommand{\R}{\mathbb{R}}
\newcommand{\N}{\mathbb{N}}

\newcommand{\Z}{\mathbb{Z}}

\theoremstyle{definition}
\newtheorem{theorem}{Theorem}

\newtheorem{proposition}{Proposition}

\newtheorem{assumption}{Assumption}
\newtheorem{lemma}{Lemma}

\title{ Detecting   Periodicity of a  General Stationary Time Series via AR(2)-Model Fitting}
\author{Jens-Peter Kreiss \footnote{Department of Mathematics, Technische Universit\"at Braunschweig, Germany, j.kreiss@tu-braunschweig.de} ,
Panagiotis Maouris 
\footnote{Department of Mathematics and Statistics, University of Cyprus, Cyprus, maouris.panagiotis@ucy.ac.cy} \ and Efstathios Paparoditis 
\footnote{Cyprus Academy of Sciences, Letters, and Arts, Cyprus, stathisp@ucy.ac.cy}}

\date{}

\begin{document}

\maketitle

\begin{abstract}
Estimating the  periodicity of a  stationary time series via fitting a second order stationary autoregressive (AR(2)) model has been initiated  by the  seminal paper of \citeasnoun{yule_method_1927}.
We  investigate properties of this procedure when applied to a   general stationary processes possessing a spectral density with a dominant peak at some frequency $\lambda_0\in(0,\pi)$. 
We  show that if the peak of the spectral density  is sharp enough (in a way to be specified)   then  the AR(2) model, which best (in  mean square sense)  approximates  the underlying  process, correctly  identifies    the  frequency $\lambda_0$.
To investigate  consistency properties of the AR(2)  based estimator  of   $\lambda_0$,    a  near to pole   framework is adopted.  Triangular arrays of  stationary stochastic processes  are considered that possess a spectral density   the   peak of which at $\lambda_0$   becomes more  
pronounced    as  the sample size $n$ of the observed time series  increases to infinity.  It is  shown  in this set up, that  the AR(2) based estimator achieves a  rate of convergence  which is larger than the parametric $n^{-1/2}$  rate and which can be  arbitrarily  close to  $ n^{-2/3}$, the best rate  that can be achieved by this estimator. 
\end{abstract}

\renewcommand{\baselinestretch}{1.2}
\small
\normalsize

\section{Introduction}

Estimating the periodicity  of a time series is a widely  discussed problem in the time series literature; see \citeasnoun{quinn_estimation_2001} for an overview. 
One well studied  case is that where the random variables $X_t$ are generated by a  nonstationary in the mean model  
that  allows for a   periodic behavior, that is,  
\begin{equation} \label{eq.meanModel-1}
X_t = \beta_1 \cos(\lambda_0t) + \beta_2 \sin(\lambda_0t) + U_t.
\end{equation}
Here $ \beta_1$, $\beta_2$  and $\lambda_0 \in(0,\pi)$ are unknown constants and $\{U_t,t\in\Z\}$ is an unobservable stationary process; see \citeasnoun{hannan_estimation_1973} and  \citeasnoun{quinn_estimation_2001}. 

Unlike the process (\ref{eq.meanModel-1}),  the focus in this paper   is  on cyclic behavior  for  weakly stationary processes $\{X_t,t\in\Z\}$. In this context a cyclic behavior   typically manifests itself through  a pronounced peak of  the  
spectral density $f$ of $ \{X_t,t\in\Z\}$  at   frequency $\lambda_0$. Recall that such a  peak  corresponds to a cycle  with period $P=2\pi/\lambda_0$. Within the class of stationary processes, 
 some different models  have been considered in the literature that are able  to describe a cyclic behavior.
A simple cosinusoidal model  is given by
\begin{equation}  \label{eq.cosin}
X_t = B_1 \cos(\lambda_0t) + B_2 \sin(\lambda_0t).
\end{equation}
Although model (\ref{eq.cosin})  looks at a first glance similar to (\ref{eq.meanModel-1}), it  differs by the fact that  the constants  $\beta_1$ and $\beta_2$ in (\ref{eq.meanModel-1})  are replaced by  random variables  $B_1$ and $B_2$. The latter  are  assumed to have mean zero, to be uncorrelated and to have  the same variance.  Model (\ref{eq.cosin})  possesses a spectral distribution function which has a jump at frequency $\lambda_0$. 

 In the time series literature,  models allowing for a pole of the spectral density at frequency $\lambda_0$ have attracted more interest. Such models can  be described by the  class of spectral densities, 
\begin{equation} \label{eq.modelGegenbauer}
f(\lambda) = \frac{1}{|1-e^{-i(\lambda-\lambda_0)}|^{2d} |1-e^{-i(\lambda+\lambda_0)}|^{2d} } f_U(\lambda), \ \ \lambda \in \R.
\end{equation}
Here,  $f_U$ refers to  the spectral density of a stationary process  $\{U_t,t\in\Z\}$ and $d$ is  a  so-called memory parameter; see among others, \citeasnoun{giraitis_gaussian_2001} and \citeasnoun{hidalgo_estimation_2004}.  Processes  possessing the spectral density  (\ref{eq.modelGegenbauer})  are said to follow a  generalized  Gegenbauer model due the appearances of the so-called Gegenbauer polynomial, $ G(z) = 1 -2z\cos(\lambda_0) + z^2 = (1-e^{-i\lambda_0}z)(1-e^{+i\lambda_0}z)$ in the expression for $f$; see \citeasnoun{gray_generalized_1989}. Notice that for $0<d<1$  and $\lambda_0 \in (0,\pi)$ or for $0<d<\frac{1}{2}$  and $\lambda_0 \in \{0,\pi\}$,  the process is stationary. \citeasnoun{giraitis_gaussian_2001}  and \citeasnoun{hidalgo_estimation_2004} considered estimation of  $\lambda_0$ using the periodogram $ I_n(\lambda) = (2\pi n)^{-1}|\sum_{t=1}^n X_t \exp\{-i t \lambda \}|^2$. In particular, they investigate $\widetilde{\lambda}_0 =\mbox{argmax}_{1\leq j\leq N}I_n(\lambda_{j,n}) $, where $ \lambda_{j,n}=2\pi j/n$,  $ j=1,2, \ldots, N=\lfloor n/2 \rfloor$ denotes the set of positive  Fourier frequencies.

The  class of  stationary autoregressive  processes also has  been  used to model periodicity. Starting point was  the seminal paper by \citeasnoun{yule_method_1927},  where a second order autoregression 
\begin{equation}\label{eq.modelAR2}
X_t = a_1 X_{t-1} + a_2 X_{t-2} + \epsilon_t,
\end{equation}
has been fitted in order  to estimate the periodicity of the time series  of sunspot-numbers; see  \citeasnoun{yule_method_1927}.  In  expression (\ref{eq.modelAR2}),  
$ \{\epsilon_t, t\in \Z\}$ has been  assumed to be a  zero mean uncorrelated (white noise) process with variance $0<\sigma_\epsilon^2<\infty$ and the parameters $ a_1, a_2 $ are such that  the stationarity and causality condition  $ a(z)=1-a_1z-a_2z^2 \neq 0$ for all $ |z|\leq 1$,   is fulfilled.  Recall that 
(\ref{eq.modelAR2}) possesses the  spectral density 
\begin{equation} \label{eq.Spec-AR2}
 f(\lambda) = \frac{\sigma^2_e}{2\pi} \Big|1-a_1e^{-i\lambda} -a_2e^{-i 2\lambda}\Big| ^{-2},
 \end{equation}
and if the roots of $a(z)$ are not real, then $a_2<0$ and  $f(\lambda)$  has   a peak (maximum) at frequency 
\begin{equation} \label{eq.lambdaAR2}
\lambda_{0}={\rm arccos}\Big(\frac{a_1(1-a_2)}{-4a_2}\Big).
\end{equation} 
If one projects the root (with positive imaginary part) of the polynomial $1-a_1z-a_2z^2$, $a_2<0$, 
on the unit circle (along the diameter), then this point relates to frequency   
\begin{equation} 
\widetilde \lambda_{0}={\rm arccos}\Big(\frac{a_1}{2\sqrt{-a_2}}\Big).
\end{equation} 
$\widetilde \lambda_{0}$ is referred in the literature as the peak frequency of the spectral density $f$ 
given in \eqref{eq.Spec-AR2}.  This
value converges to the maximal value of $f$ given in \eqref{eq.lambdaAR2} if the roots of the polynomial converge to the
unit circle.
Plugging in (\ref{eq.lambdaAR2}) estimators of the parameters  $a_1$ and $a_2$  leads to a (moment) estimator of  $\lambda_0$, since $ a_1$ and $a_2$ and consequently  $ \lambda_0$,  can be expressed as functions of the lag one and lag two autocorrelations through  solving  the  well-known Yule-Walker equations; see \citeasnoun{brockwell_time_1991}, Chapter 8.  In contrast to the  jump behavior in (\ref{eq.cosin}) or the pole behavior in  (\ref{eq.modelGegenbauer}), the spectral density  of the AR(2) model  (\ref{eq.modelAR2}) is analytic even at the peak $\lambda_0$.  

The aim of this paper is to  thoroughly investigate the capability  and the consistency  properties of such an AR(2) model  
fit for detecting the  periodicity $\lambda_0$  when the underlying process $\{X_t,t\in\Z\}$  is not necessarily an  AR(2) process (see equation  (\ref{eq.modelAR2})),  but  
a  general stationary process. 
Inspired by the set up   in \citeasnoun{dahlhaus_1988}, we focus on the situation where  the  stationary process considered possesses  a sharp  peak at  some frequency $\lambda_0\in(0,\pi)$.  Note that the standard asymptotic set up is not capable to properly resemble  the inference problem we have in mind. To elaborate,  it is well known that if an AR(2) model is fitted to a  general 
stationary time series by means of minimizing the mean square  error, then  the Yule-Walker parameter estimators converge to their theoretical counterparts at $n^{-1/2}$ rate and this rate leads, under certain conditions,  to a $ n^{-1/2}$ rate  also for the corresponding 
estimator of $\lambda_0$.  However, this set up does not allow to investigate  how and to what extent  the sharpness of the peak  affects the capability of the AR(2) model to detect periodicity and   how this sharpness 
impacts   the corresponding inference problem and the rate of convergence of the estimator.  Recall that   the situation  we consider is that  where    the spectral density of the underlying stationary process  has 
  a  strong and dominant peak at  some (unknown) frequency $\lambda_0$. 
  
   To develop a  more appropriate mathematical framework for our  investigations, 
we  consider  a class of stationary 
 processes  possessing  a   spectral  density  that can be expressed as   
\begin{equation} \label{eq.f-delta}
f_\delta(\lambda) = \frac{1}{|1-(1-\delta)e^{-i(\lambda-\lambda_0)}|^{2} |1-(1-\delta)e^{-i(\lambda+\lambda_0)}|^{2} } f_U(\lambda),
\end{equation}
where $\delta$ is a parameter  taking  values in the interval   $ (0,1)$.   The  class of  spectral densities given in   (\ref{eq.f-delta}) generalizes   in some sense  the spectral densities  considered so far  
and has a  number of interesting features. Note  that  $f_\delta$ has a peak at frequency $\lambda_0$ which becomes the more  pronounced the closer   $\delta$  is to zero.  In fact  the parameter $\delta$ allows to   control  within the stationary process class, the 
strength of the peak   of $ f_\delta$ at  $\lambda_0$.  For $\delta $ tending to  zero, the behavior of $f_\delta$ at $\lambda_0$  gets closer to   that of  a pole, that is $f_\delta$ 
resembles  in the limit,  the behavior   of   the spectral density (\ref{eq.modelGegenbauer}). Apart from  the (long) memory properties described through the parameter $d$ and which are not the focus of the current paper,  the pole behavior of 
(\ref{eq.modelGegenbauer}) at $\lambda_0$  can therefore    be considered  as  a limiting case    for  $\delta\rightarrow 0$ of  the  class of  spectral densities given in  (\ref{eq.f-delta}).  On the other hand,  for $ \{U_t,t\in\Z\}$ a white noise process and for $ a_1=2(1-\delta)\cos(\lambda_0)$ and $ a_2=-(1-\delta)^2$, the spectral density $  f_\delta $ coincides  for any fixed $\delta\in(0,1)$, with that of a pure AR(2) model as given in (\ref{eq.Spec-AR2}).  Notice that  if  the autoregressive polynomial $ a(z)=1-a_1z - a_2z^2$ has two complex roots   $z_{1,2}=x\pm i y $, then  $ a_1=2(1-\delta)\cos(\lambda_0)$ and $ a_2=-(1-\delta)^2$, where $ \delta =1-\sqrt{x^2+y^2}$ and $ \lambda_0= \tan^{-1}(y/x)$.   As we will see, the framework defined by the spectral density class (\ref{eq.f-delta}),   also enables us to investigate the effects of the sharpness of the peak on the rate of convergence of the AR(2)-based sample estimator of $\lambda_0$. To achieve this  we will  allow   $ \delta$  to depend on the length $n$ of the time series at hand and  will let   $\delta$  getting  closer to zero, that is, we will make   the peak of the spectral density sharper,  as $ n$ increases to infinity.

A special case of   (\ref{eq.f-delta}), in which $ \delta=0$ and $f_U$ is  the spectral density of a linear process driven by i.i.d. innovations,  has attracted  more interest in the econometric time series literature.  In this case, 
the  corresponding process is a non-stationary, linear process which   possesses   two (conjugate)  complex  roots on the unit disc;  see \citeasnoun{AhtolaTiao1987}, \citeasnoun{ChanWei1988} and \citeasnoun{Tanaka2008AnalysisOM}. 

For   general stationary processes possessing  the  spectral density (\ref{eq.f-delta}) and to the best of our knowledge, the effectiveness  of the AR(2) model to  correctly identify the  periodicity and  of the  corresponding sample estimator  to   consistently estimate the peak frequency $\lambda_0$, are  largely unknown in the literature.  More specifically,   it is not known if  and under what circumstances  
an AR(2) fit to the spectral density class (\ref{eq.f-delta})  correctly identifies the frequency  $\lambda_0$.  An affirmative   answer to   this question turns out to be  interesting  especially    when taking   into account the fact that the  effectiveness  of such an  AR(2) fit  has been largely questioned in the literature. In particular,  \citeasnoun{quinn_estimation_2001} have shown, that in the context of the non stationary in the mean model (\ref{eq.meanModel-1}), a pure   AR(2) fit,  succeeds  in  correctly identifying the frequency  $\lambda_0$,  if and only if  the process $ \{U_t, t\in\Z\}$ in (\ref{eq.meanModel-1})  has variance zero. 

Apart from this, also the   consistency  properties  of  an   AR(2) based estimator of $\lambda_0$ for general stationary processes   possessing a spectral density like in  (\ref{eq.f-delta})  are not known.  The setting of an early work  by   \citeasnoun{newton_method_1983} where the authors  consider  fitting  a  general $p$th order autoregression for peak estimation,  is   far too restrictive for our purposes. In particular, the authors  assume    that  $\{X_t,t\in\Z\}$ follows  a finite or infinite order, (linear) Gaussian  autoregression  and the  $n^{-1/2}$-consistency  rate of the corresponding  frequency estimator has been  derived  under  the assumption that the order $p$ of the underlying process is  finite and known. This rate has been weakened to  $p^{3/2}n^{-1/2}$, in  case $\{X_t,t\in\Z\}$ is  an  infinite order autoregressive process  and  the order $p$ of the autoregressive model    fitted,  increases  to infinity with  the  sample size $n$.  As  already mentioned,  we do not  aim  to  impose  structural assumptions  on  the  stationary process $\{U_t,t\in\Z\}$ involved. Moreover,  in the context considered in this paper,    the  AR(2) model   is solely used   as a vehicle in order to  identify and infer  periodicity, i.e., to identify and estimate $\lambda_0$  and not   as a procedure to parametrize the  dependence structure of the underlying process.

In Section 2  we investigate within the process class (\ref{eq.f-delta}),  the effectiveness of the AR(2) model  to  correctly identify the frequency $\lambda_0$ of interest.  For this  we  consider the case where the peak of the process at $\lambda_0$  is sharp enough  and this is  achieved  by    allowing  the parameter  $ \delta$ to converge to zero. We   show that  if  $\delta$ is  close enough to zero and despite the properties that  the spectral density  $f_U$ may have,  
the frequency  $\lambda_{0,\delta} $ where the spectral density of the best (in mean square sense) approximating   AR(2) model  has its peak, is   close enough  to  $\lambda_0$. In fact  the difference between $ \lambda_{0,\delta}$ and $\lambda_0$ can be arbitrary small and  in the limit this discrepancy  converges to  zero for  $\delta\rightarrow 0$.   In Section 3 we investigate the consistency   properties of a sample  estimator of $ \lambda_{0,\delta}$ obtained by  fitting an AR(2) model to a time series $ X_1, X_2, \ldots, X_n$ stemming from a stationary process with  spectral density in the class (\ref{eq.f-delta}).  The estimators of the autoregressive parameters are obtained   via solving the system of Yule-Walker equations. 
To develop a theory capable  to investigate  properties of this  estimator,  we adopt a ``close to pole'' asymptotic framework. For this,  
we  allow   the parameter  $ \delta$  to  depend on the sample size $ n$   such that    $\delta=\delta(n) \rightarrow 0$ as $ n \rightarrow \infty$. That is,  we allow for  the  peak of   $ f_{\delta(n)}$  at $\lambda_0$ to become more and more pronounced the larger   the sample size   gets.  By properly controlling  the rate at  which $\delta(n) $  approaches zero with respect to $n$, we show that the AR(2)-based estimator  is consistent  and that the  rate of convergence  of this estimator  can be arbitrarily close to  $ n^{-2/3}$.

Before proceeding with our investigations, we fix some notation. For two sequences $(a_n) $ and $ (b_n)$, we use the notation $ a_n\asymp O(b_n)$ if  constants  $k_1>0$ and $k_2>0 $ exist such that $ k_1 b_n \leq a_n \leq k_2 b_n$.  



\section{Effectiveness of the AR(2) Model in Detecting Periodicity}

We start with the  following assumption which describes our requirements on the underlying process $ \{X_{t,\delta}, t \in\Z\}$. 
 
\begin{assumption} \label{as.1}
 For $\delta\in(0,1)$ and $\lambda_0\in(0,\pi)$,   $\{X_{t,\delta}, t\in \Z\}$  is  a  zero mean stationary stochastic process which possesses  a   spectral density $f_\delta(\lambda)$, $\lambda \in \R$, as  given in  (\ref{eq.f-delta}). 
$ f_U $  is the spectral density of a stationary process $\{U_t,t\in\Z\}$, which is continuous everywhere in $[0,\pi]$ and 
therefore also bounded, and satisfies $ \inf_{\lambda\in[0,\pi]}f_U(\lambda) >0$. 
 \end{assumption}
 
 \vspace*{0.2cm}
 
 \noindent{\bf Remark:}
 \begin{enumerate}
 \item[(i)] \ 
 Assumption~\ref{as.1}  is fulfilled  if $ X_{t,\delta}$  can be expressed as 
 \begin{equation}
 (1 - (1-\delta)e^{-i\lambda_0}L)(1-(1-\delta)e^{i\lambda_0}L) X_{t,\delta} =  U_t,
 \end{equation}  
 where $ L$ is the shift operator defined as $ L^kX_{t,\delta} = X_{t-k,\delta}$ for $k\in \Z$ and $ \{U_t,t\in\Z\}$ is a  zero mean stationary process possessing the  spectral density $ f_U$. 
  Notice that  
   $ |(1-\delta)e^{\pm i\lambda_0}|<1$, that is, the process $\{ X_{t,\delta},t\in\Z\}$  is stationary  and  causal with respect to the process $\{U_t,t\in \Z\}$.  $X_{t,\delta}$ 
    also can be  expressed  as $X_{t,\delta}=A_\delta(L)^{-1} U_t$, where  
  $ A_\delta(L)=(1-2(1-\delta) \cos(\lambda_0)L + (1-\delta)^2 L^2)$.
 \item[(ii)] \ The assumptions imposed on the spectral density $f_U$ imply that $\{U_t,t\in\Z\}$ has a so-called autoregressive representation, that is, a sequence of coefficients $ \{c:j\in {\mathbb N}\}$ exists with $ \sum_{j=1}^\infty |c_j|<\infty$ and $ c(z)=1-\sum_{j=1}^\infty c_j z^j\neq 0$ for $ |z|\leq 1$, such that $ U_t=\sum_{j=1}^\infty c_j U_{t-j} + e_t$.  Here $ \{e_t,t\in\Z\}$ is a white noise process, that is, ${\rm  E}(e_t)=0$ and $ \cov  (e_t,e_s)=0$ for $ t\neq s$ and ${\rm Var}(e_t)=\sigma^2_e$.  Hence 
 $ A_\delta(L)C^{-1}(L) X_{t,\delta} = e_t$, where $ C^{-1}(z)$ denotes the inverse power series of $C(z)=1+\sum_{j=1}^\infty c_jz^j$, $ |z| <1$. 
 \end{enumerate}

 \noindent For the stationary process $ \{U_t,t\in\Z\}$ we additionally assume that its second and fourth order cumulants 
 satisfy the following summability conditions.
 
\begin{assumption} \label{as.2}
For the autocovariance function of $ \gamma_U(h) = {\rm Cov}(U_{t}, U_{t+h})$, $h\in\Z$,  and the cumulant function   
$ \textrm{cum}(U_{0}, U_{h_1,}, U_{h_2,}, U_{h_3} )$, $ h_1,h_2,h_3\in\Z$,  of the stationary process $ \{U_t,t\in\Z\}$, 
 it holds true that 
\begin{enumerate}
    \item[(i)] \  $\sum_{h \in \mathbb{Z}}  |h||\gamma_{U}(h)|  < \infty$.
    \item[(ii)] \  $\sum_{h_1, h_2, h_3 = -\infty}^{\infty} (1+|h_1|+|h_2|+|h_3|) |\textrm{cum}(U_{0}, U_{h_1,}, U_{h_2,}, U_{h_3} )|  
    < \infty$\\
\end{enumerate} 
\end{assumption}

 \noindent As already mentioned, the spectral density $f_\delta$ has a peak at frequency $\lambda_0$ which is the sharper the closer $\delta $ is  to zero.   Let  $ \gamma_\delta(h)={\rm Cov}(X_{\delta,t},X_{\delta,t+h})$ and $ \rho_\delta(h) = \gamma_\delta(h)/ \gamma_\delta(0) $, $ h \in \Z$, be the autocovariance and the autocorrelation function, respectively,  of $\{X_{t,\delta},t\in\Z\}$.
 The following lemma deals with  the behavior of  these second order  functions as well as with some  basic properties of the model class (\ref{eq.f-delta}).

 \begin{lemma}  
 \label{le.1}
Under Assumption~\ref{as.1}, it holds 
\begin{enumerate}
 \item[(i)] \ $    \lim_{\delta\rightarrow 0}f_\delta(\lambda_0)=+\infty  \ \mbox{and} \ f_\delta (\lambda) \asymp {\mathcal O}(\delta^{-2})  \ \ \mbox{as} \ \ \lambda\rightarrow \lambda_0$
  \item[(ii)] \ $  \lim_{\delta\rightarrow 0}\gamma_{\delta}(0) =   + \infty  \ \   \mbox{and} \ \ \gamma_{\delta}(0)  \asymp {\mathcal O}(\delta^{-1})  \  \ \mbox{as} \ \ \delta\rightarrow 0$
   \item[(iii)] \  $ \rho_{\delta}(h)\ = \cos(\lambda_0 h) + {\mathcal O}(\delta)$ as $ \delta \rightarrow 0$ for any  $   h\in \mathbb{N}_0 $
    \item[(iv)] \  $ X_{t,\delta}=\sum_{j=1}^\infty b_{j,\delta} U_{t-j} + U_t,  \ \mbox{where for} \  K=\lfloor j/2\rfloor -1$,
    \[    b_{j,\delta} = \left\{ \begin{array}{lll} 
      (1-\delta)2\cos(\lambda_0) & & j=1 \\
       (1-\delta)^{j} \big( 1+2\cos(j\lambda_0) - \frac{\displaystyle 2\sin((K+1)\lambda_0)\cos(K\lambda_0)}{\displaystyle \sin(\lambda_0)}\big)& & j>1 \ \mbox{even}\\
         (1-\delta)^{j} \big( -1+2\cos(j\lambda_0) - \frac{\displaystyle 2\cos((K+1)\lambda_0)\sin(K\lambda_0)}{\displaystyle \sin(\lambda_0)} +\frac{\displaystyle \sin(4K+3)\lambda_0/2)}{\displaystyle \sin(\lambda_0/2)} \big)& & j>1 \ \mbox{odd}
        \end{array} \right.  \]
     \mbox{In  particular,}  $ |b_{j,\delta}|   \leq C_{\lambda_0}(1-\delta)^j \  \mbox{for all} \   j \in \N \
     \mbox{where the constant}  \ C_{\lambda_0}>0,\\
       \  \mbox{only depends on} \ \lambda_0$.
\end{enumerate}
\end{lemma}

 \noindent As the above lemma shows the spectral density of the process produces a pole-similar behavior at $\lambda=\lambda_0$ as $ \delta $ approaches zero  while   the  variance of  $\{X_{t,\delta}, t \in \Z\}$ diverges to infinity.
At the same time and as $ \delta \rightarrow0$, the behavior of the autocorrelation function   gets closer to that of a  $\cos(\cdot)$ function. As a matter of  fact, the limiting autocorrelation function coincides 
 with the  autocorrelation function of  the  cosinusoidal model (\ref{eq.cosin}). 
 
\noindent The following  result  answers the question about  the effectiveness   of the best (in means square sense) fitting  AR(2) model  to detect the peak frequency  $\lambda_0$ of interest.  

\begin{lemma} \label{le.2}
Suppose that the  AR(2) model  $ X_{t,\delta} = a_{1,\delta} X_{t-1,\delta} + a_{2,\delta} X_{t-2,\delta} + v_t$ is fitted to the process $\{X_{t,\delta}, t\in \Z\}$, where   the unique  parameters
$ (a_{1,\delta}, a_{2,\delta})$ are those obtained as $ {\rm argmin}_{c_1,c_2} {\rm E}(X_{t,\delta}-c_1X_{t-1,\delta} -c_2X_{t-2,\delta})^2$. 
If  $ a_{2,\delta}<0$  define  $ \lambda_{0,\delta} = \arccos\big(a_{1,\delta}(1-a_{2,\delta})/(-4a_{2,\delta})\big)$.
Then,   
\begin{enumerate}
\item[(i)]  \  $ \lim_{\delta \rightarrow 0}\  (a_{1,\delta}, a_{2,\delta}) = (2\cos(\lambda_0) , -1)$
\item[(ii)]  \  $  \lim_{\delta\rightarrow 0} \lambda_{0,\delta} = \lambda_0$.
\end{enumerate}
\end{lemma}
\vspace*{0.25cm}

\noindent Recall that the parameters $ a_{1,\delta}$ and $a_{2,\delta}$  of the best AR(2)-fit are functions of the first and second order autocorrelations and they are  given by
\begin{equation} \label{eq.YW-Est}
(a_{1,\delta}, a_{2,\delta}) = \frac{1}{1-\rho^2_{\delta}(1)}
\big( \rho_{\delta}(1)(1-\rho_{\delta}(2)),  \   \rho_{\delta}(2)-\rho^2_{\delta}(1)\big).
\end{equation}
$a_{1,\delta}$ and $ a_{2,\delta}$  are obtained  by solving the corresponding two-dimensional system of 
Yule-Walker equations; see \citeasnoun{brockwell_time_1991}, Chapter 6.

Lemma~\ref{le.2} justifies, under certain conditions,  the use of an   AR(2) model in order  to  detect the dominant  periodicity of a time series at hand.
Recall that the closer is   $\delta$    to zero 
the more pronounced is the peak of the spectral density $ f_\delta$ at frequency  $\lambda_0$.  Therefore,  Lemma~\ref{le.2} states that   if the peak of the spectral density is strong enough, the frequency  $\lambda_{0,\delta}$ where the  best (in mean square sense) fitting AR(2) model has its peak  is close enough to the frequency $\lambda_0$.

\section{Consistency of the AR(2)-Based  Estimator of Periodicity}

Suppose a time series  $ X_{1,\delta}, X_{2,\delta}, \ldots, X_{n,\delta}$  stemming from  the  process  $\{X_{t,\delta}, t\in \Z\}$  is observed. Since in  our consistency considerations, we do not want to impose any  structural  assumptions on the dependence structure of the process $\{U_t,t\in\Z\}$ beyond the summability conditions stated in Assumption~\ref{as.2}, it  seems  more convenient to switch to the frequency domain. Towards  this goal,  the sample estimator
\[ \widehat{\gamma}_\delta(h) = \frac{2\pi}{n}\sum_{j\in {\mathcal G}(n)} I_{n,\delta}(\lambda_{j,n}) \cos(\lambda_{j,n} h),\ \ \  {\mathcal G}(n)=\{-N,\ldots, -1, 1, \ldots, N\},\]
of $\gamma(h)$ for $ |h| <n$,  is used, where $ N=\lfloor n/2\rfloor$ and  $ I_{n,\delta}$ is the periodogram of the time series   $ X_{1,\delta}, X_{2,\delta}, \ldots, X_{n,\delta}$ given by 
\[ I_{n,\delta}(\lambda_{j,n}) = \frac{1}{2\pi n} \big| \sum_{t=1}^n X_{t,\delta} e^{-i\lambda_{j,n} t}\big|^2,\] 
with $ \lambda_{j,n}=2\pi j/n$, $ j \in {\mathcal G}(n)$. Notice that  the estimator $ \widehat{\gamma}(h)$ is closely related to the commonly used sample autocovariance 
$ \widetilde{\gamma}_\delta(h) =n^{-1}\sum_{t=1}^{n-|h|} X_{t,\delta}X_{t+|h|,\delta}$, $ |h|<n$. 

\begin{proposition}\label{prop.1} Suppose that Assumption~\ref{as.1} is satisfied and that $\sum_{h \in \mathbb{Z}} |\gamma_{U}(h)|  < \infty$.  Then, for any $ |h|<n$,  it holds true that 
\begin{equation} \label{eq.Est-Gammah}
\widehat{\gamma}_\delta(h) = \widetilde{\gamma}_\delta(h) + O_{P}(n^{-1}\delta^{-2}).
\end{equation}
\end{proposition}
Clearly 
$ |  \widehat{\gamma}_\delta (h) - \widetilde{\gamma}_\delta (h)| \stackrel{P}{\rightarrow} 0$   if   $ n\delta^2 \rightarrow \infty$, as  $ n\rightarrow \infty$,    It is worth mentioning here that as an  inspection of the proof of Proposition~\ref{prop.1} shows,  for  any fixed $ \delta\in(0,1)$, 
the difference between the two estimators,  which are  denote in this case by $ \widehat{\gamma}(h)$ and $ \widetilde{\gamma}(h)$, respectively,  satisfies 
$$\widehat{\gamma}(h) = \widetilde{\gamma}(h) + O_{P}(n^{-1}).$$

Let $ (\widehat{a}_{1,\delta}, \widehat{a}_{2,\delta})$ be the well-known Yule-Walker estimator of $ (a_{1,\delta}, a_{2,\delta})$, i.e., the estimator obtained by replacing $ \rho_\delta(1)$ and $ \rho_\delta(2)$ in (\ref{eq.YW-Est}) by the corresponding sample estimators $ \widehat{\rho}_\delta(h) =\widehat{\gamma}_\delta(h)/\widehat{\gamma}_\delta (0)$.
For $\widehat{\alpha}_{2,\delta}<0$, the corresponding AR(2) based 
estimator of $\lambda_{0,\delta}$  is then given  by   
\begin{equation} \label{eq.lambda-estimator}
\widehat{\lambda}_{0,\delta} = \arccos\Big(\frac{\widehat{a}_{1,\delta}(1-\widehat{a}_{2,\delta})}{-4\widehat{a}_{2,\delta}}\Big).
\end{equation}


\noindent To investigate the properties of the estimator  (\ref{eq.lambda-estimator}), we first establish the following lemma which refers to certain summability properties of the second and fourth order moments
of the process  $ \{X_{t,\delta}, t\in \Z\}$.  

\begin{lemma} \label{le.Sum} 
  Under Assumption~\ref{as.1} and  Assumption~\ref{as.2}   it holds true that 
\begin{enumerate}
  \item[(i)] \  $\sum_{h \in \mathbb{Z}}  |\gamma_{\delta}(h)|  = {\mathcal O}(\delta^{-2})$\\
       \item[(ii)] \  $\sum_{h \in \mathbb{Z}}  |h||\gamma_{\delta}(h)| = \mathcal{O}(\delta^{-2})$\\
    \item[(iii)]\  $\sum_{h_1, h_2, h_3 = -\infty}^{\infty} (1+|h_1|+|h_2|+|h_3|) |\textrm{cum}(X_{0,\delta}, X_{h_1,\delta}, X_{h_2,\delta}, X_{h_3,\delta} )|  = \mathcal{O}(\delta^{-4})$ \\
\end{enumerate}
\end{lemma}

\noindent The next  lemma deals with  the moment structure of the periodogram $ I_{n,\delta}$  of a time series $X_{1,\delta}, X_{2,\delta}, \ldots, X_{n,\delta} $ stemming from the process $ \{ X_{t,\delta}, t \in \Z\}$.

\begin{lemma} \label{le.Per}
  Under Assumptions~\ref{as.1} and \ref{as.2}   it holds true that  
  \begin{align}
  &\text{(i)}\  E[I_{n,\delta}(\lambda_{j,n})] =  f_{\delta}(\lambda_{j,n}) + \mathcal{O}(\frac{1}{n\delta^2})\\
    &\text{(ii)}\  \cov [I_{n,\delta}(\lambda_{k,n}), I_{n,\delta}(\lambda_{m,n})] =  \left\{ \begin{array}{lll}
f^2_{\delta}(\lambda_{k,n}) + \mathcal{O}\Big(\frac{\displaystyle 1}{\displaystyle |A_{\delta}(e^{-i\lambda_{k,n}})|^2 n\delta^2}\Big) + \mathcal{O}\Big(\frac{\displaystyle 1}{\displaystyle n^2\delta^4}\Big) & & \\ 
\ \hspace*{5.5cm} \mbox{if} \ \ \   0<\lambda_{k,n}=\lambda_{m,n}<\pi & &  \\
& & \\
\frac{\displaystyle 2\pi}{\displaystyle n|A_{\delta}( e^{-i\lambda_{k,n}})|^2 |A_{\delta}(e^{-i \lambda_{m,n}})|^{2}}  f_{4,U}(\lambda_{k,n}, -\lambda_{m,n}, \lambda_{m,n})    & & \\
  \ \ \ \ + \ \mathcal{O}\big(\frac{\displaystyle 1}{\displaystyle n^2\delta^4}\big) \hspace*{3.7cm}  \mbox{if} \ \ \ 0<\lambda_{k,n}\neq\lambda_{m,n}<\pi & &
\end{array} 
\right.
\vspace{-0.1cm}
\end{align}
\hspace*{0.75cm}where  
\begin{align*}
f_{4,U}(\lambda_{k,n}, -\lambda_{m,n}, \lambda_{m,n})
= (2 \pi)^{-3} \sum_{h_1,h_2,h_3 \in \mathbb{Z}} \textrm{cum} (U_{0}, U_{h_1},U_{h_2}, U_{h_3}) e^{-i(h_1\lambda_{k,n}-h_2\lambda_{m,n}+h_3\lambda_{m,n})} 
\end{align*}
is the fourth order cumulant spectral density of the stochastic  process  $ \{U_t,t\in\Z\}$, which because of  Assumption~\ref{as.2}
exists.
\end{lemma}

\noindent As mentioned in the Introduction, in order to develop   a meaningful asymptotic theory which is capable to investigate the consistency  properties of the estimator (\ref{eq.lambda-estimator}),  we adopt  a  ``close to pole'' asymptotic framework. For this,  
we   allow     $ \delta$ to  depend  on the sample size $ n$ such that  $\delta_n:=\delta(n) $ converges to zero as  $ n$ increases to infinity at some  controlled rate, which is  stated  in  the following assumption.

\vspace*{0.2cm}

\begin{assumption} \label{as.3}  $ \delta_n \sim n^{-\alpha} $ for some  $ \alpha\in(0,1/3)$. 
\end{assumption}

\vspace*{0.2cm}
 
\noindent Consider now an array of stochastic processes,  denoted by $ \{X_{t,\delta_n}, t \in \Z\}_{n\in\N}$,  such that for each $ n\in \N$, $ \{X_{t,\delta_n}, t\in\Z\}$  possesses the spectral density $ f_{\delta_n}$ 
given in (\ref{eq.f-delta}). For   $ X_{1,\delta_n}, X_{2,\delta_n}, \ldots, X_{n,\delta_n}$ the  time series available,  denote by $ \widehat{\lambda}_{0,\delta_n} $ the estimator (\ref{eq.lambda-estimator}) based on this time series. 
The following  consistency result can then be established.

\begin{lemma} \label{le.gamma0}
  Under Assumptions~\ref{as.1}, ~\ref{as.2} and  ~\ref{as.3},  it holds true, as $n\rightarrow \infty$, that 
  \begin{equation*}
      \frac{\hat{\gamma}_{\delta_n}({0})}{\gamma_{\delta_n}({0})} \xrightarrow {P} 1
  \end{equation*}
\end{lemma}

\noindent  We can now state our main result describing  the consistency and the rate properties of the parameter estimators of  the AR(2) fit as well as of the corresponding estimator of the
location of the peak.

\begin{theorem} \label{th.1}
Suppose that   Assumptions~\ref{as.1}, ~\ref{as.2} and ~\ref{as.3} are satisfied. Then,  as $ n\rightarrow\infty$,  
\begin{enumerate}

\item[(i)] \ 
$ \widehat{\alpha}_{1,\delta_n}- \alpha_{1,\delta_n} ={\mathcal O}_{P}\big(n^{-(1+\alpha)/2}\big)  $ 
and $ \widehat{\alpha}_{2,\delta_n}- \alpha_{2,\delta_n}={\mathcal O}_{P}\big(n^{-(1+\alpha)/2}\big) $.
\item[(ii)] \
$ \widehat{\lambda}_{0, \delta_n} = \lambda_{0,\delta_n} \ +\  {\mathcal O}_P  (n^{-(1+\alpha)/2})$.
\end{enumerate}
\end{theorem}

\noindent {\bf Remark.}  
\noindent Recall that $ \lambda_{0,\delta_n}$ is the location of the peak of the spectral density of the AR(2) model which 
best (in mean square sense) approximates  the second order structure of the underlying process $\{X_{t,\delta_n},t\in\Z\}$;  
see Lemma~\ref{le.2}.  By Theorem~\ref{th.1},  the sample estimator  $\widehat{\lambda}_{0,\delta_n}$ of the peak 
location  converges at a rate which is  larger than the parametric rate $n^{-1/2}$ and 
this  rate is the closer to $n^{-2/3}$ the  stronger the peak of the underlying spectral density  is. \\
  
Giraitis et al. (2001) obtained for the maximizing frequency $\lambda_0$ of their objective function over a discrete set of frequencies and for 
infinite order moving average models  with ergodic innovations which satisfy some conditional moment assumption up to 
order four, a convergence rate of  $n^{-1}$.   
Hidalgo and Soulier (2004) considered the maximum of the periodogram of all Fourier frequencies as an estimator for
the spectral peak. They were able to show  for linear time series, i.e. moving averages of infinite order driven by  i.i.d. innovations, 
a convergence rate of their estimator which almost reaches $n^{-1}$.
The rate we obtained  for the AR(2) based estimator is slower.  However,  on the one hand we think that it is worth  investigating   the ability 
of the   simple AR(2)-based procedure  to  estimate the peak frequency  under rather mild assumptions,  an issue which,  to the best of our knowledge,  has not been  addressed  so far in the literature.
On the other
hand,  we conjecture that an asymptotic distribution theory may be possible for the AR(2)-based estimator of the
spectral peak, a task which so far, is open for further research. This could be an advance of the AR(2)-based estimator  taking into account 
the fact that  for the estimator  based on the maximum of the periodogram  over the Fourier frequencies, a distribution  theory is not 
possible; see  Giraitis et al. (2001), p. 997.


%
%
\section{Estimating the Periodicity of the Sunspot Series}
\begin{figure}[h]
    \centering
    \includegraphics[width=12.5cm,height=0.2\textheight]{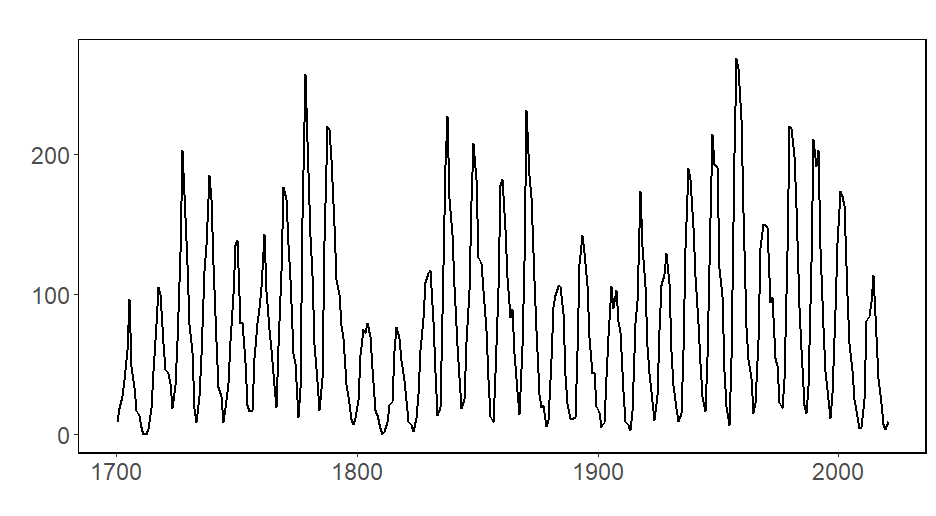}
     \begin{center}
    \begin{minipage}{11cm}
    \vspace*{-0.5cm}
    \it
    \caption{Plot of the times series of  Sunspot  Numbers   from 1700 to 2020.}
    \label{fig:Sunspot}
        \end{minipage}
    \end{center}
\end{figure}

\begin{figure}[h]
    \centering
    \includegraphics[width=12.5cm,height=0.2\textheight]{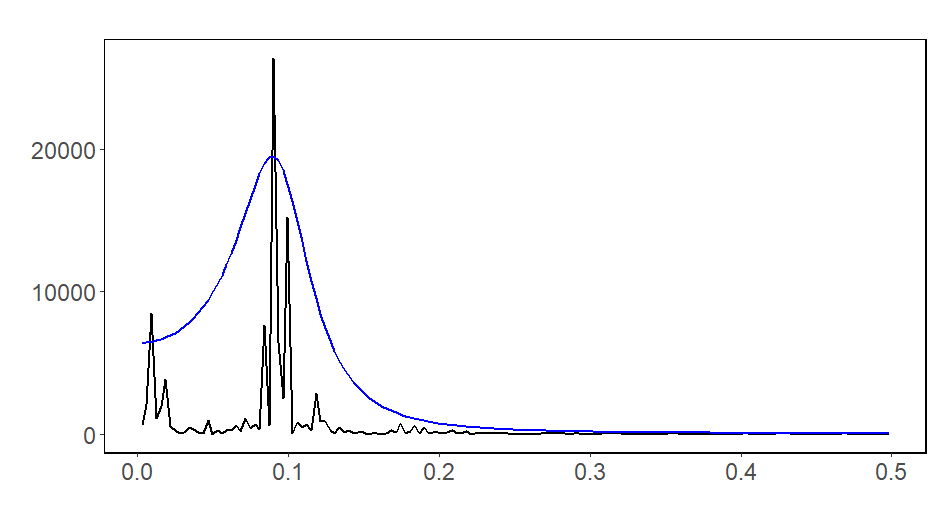}
    \begin{center}
    \begin{minipage}{11cm}
    \vspace*{-0.5cm}
    \it
    \caption{Periodogram of Sunspot Series (black solid line) and 
    the Spectral Density of the fitted AR(2) model (blue line).}
    \label{fig:fitAR2}
    \end{minipage}
    \end{center}
\end{figure}
%

In this section we apply the procedure investigated in this paper to estimate  the periodicity of  the yearly mean sunspot numbers observed from 1700 to 2020 (n = 321).  A plot of the time series 
which  is available at  
 www.sidc.be/silso/datafiles,  is shown in 
Figure~\ref{fig:Sunspot}.  The data measures the solar activity and specifically the number of sunspots visible on  sun's surface. It is calculated by counting the number of individual sunspots and sunspot groups and by  combining those counts in specific way. The sunspot numbers are used to track the solar cycle of increasing and decreasing solar activity.

 Fitting a second order autoregressive model to this time series leads to the estimates 
 $\widehat{a}_1=1.386$ and $ \widehat{a}_2=-0.691$ of the corresponding parameters.  Furthermore,  $ \widehat{\alpha}_1^2- 4\widehat{\alpha}_2 =-0.843$, which shows that the second order  characteristic polynomial of the fitted AR(2) model has two
 (conjugate) complex roots.  Based on the estimates of the AR(2) parameters we  get  using formula (\ref{eq.lambdaAR2}), the estimate  $\widehat{\lambda}_0= 0.559$  of the frequency where the spectral density of the fitted AR(2)  model  has its peak. This   corresponds to an estimated periodicity of $ \widehat{P}=11.24$ years for  the sunspot series.
 It is worth  mentioning that  estimating  the peak frequency  using the maximum of the periodogram over the Fourier frequencies, leads to the estimate $\widetilde{\lambda}_{0} = 0.568$ which  corresponds to an estimated  periodicity  of $\widetilde{P}= 11.069$ years,  which is  quite close to the estimate of the periodicity  obtained using the AR(2) fit.  Figure~\ref{fig:fitAR2}  shows the periodogram of the time series of sunspot numbers together with a plot of the spectral density of the fitted AR(2) model.  \\
%
%

\section{Technical Results and Proofs}

\noindent {\bf Proof of Lemma~\ref{le.1}:}  
(i) Using 
\begin{align*} f_\delta(\lambda)= & \frac{1}{\big(1-2(1-\delta)\cos(\lambda-\lambda_0)+(1-\delta)^2\big)} \\
& \times \frac{1}{\big(1-2(1-\delta)\cos(\lambda+\lambda_0)+ (1-\delta)^2\big)} \cdot f_U(\lambda),
\end{align*}
and  $ \cos(2x)=2\cos^2(x)-1$ we   get for $\lambda=\lambda_0$,
\[ f_\delta(\lambda_0) =\frac{1}{ \delta^2}\cdot \frac{1}{4(1-\delta)[1-\cos^2(\lambda_0)] +\delta^2}\cdot f_U(\lambda_0).\] 
(ii) Recall that 
\begin{align*}
  \gamma_{\delta}(0)\ &= 2\int^{\pi}_{0} f_{\delta}(\lambda)\ d\lambda \\
&= 2\int_{0}^\pi |1-(1-\delta)   e^{-i(\lambda\ -\lambda_0)}|^{-2} \frac{ f_{U}(\lambda)}{ 1-2(1-\delta) \cos(\lambda + \lambda_0) + (1-\delta)^2 } d\lambda
 \end{align*}
 and observe that $ c_2 >c_1 >  0$   exist such that for $\lambda \in [0,\pi ]$
 \begin{equation*}
  c_1\cdot|1-(1-\delta)e^{-i(\lambda\ -\lambda_0)}|^{-2} \leq f_{\delta}(\lambda) \leq c_2\cdot|1-(1-\delta)e^{-i(\lambda\ -\lambda_0)}|^{-2}.
 \end{equation*}
 To see this note  first that since  $ \lambda_0\in(0,\pi)$,  $$
  -1 \leq  \cos(\lambda + \lambda_0) \leq  c_{\text{max}}$$ for $ \lambda \in [0,\pi]$ 
 where  $c_{\text{max}} = \text{max}\{\cos(\lambda_0), \cos(\pi + \lambda_0)\} \in [0,1)$. Using   $ \cos(\lambda + \lambda_0)\geq -1$ we easily get 
 \begin{equation} \label{eq.B-1}
 1-2(1-\delta) \cos(\lambda + \lambda_0) + (1-\delta)^2 \leq  1+2(1-\delta) + (1-\delta)^2 = (2-\delta)^2 \leq 4.
 \end{equation}
 Furthermore, using 
 \[ 1-2(1-\delta)\cos(\lambda+\lambda_0) +(1-\delta)^2 \geq  1-2(1-\delta)c_{\max}+(1-\delta)^2,\]
 and since the function on the  right hand side of above   inequality has the  unique   minimum $1-c_{\max}^2 \in (0,1]$ for $ \delta=1-c_{\max}$,  we get  
 \begin{equation} \label{eq.B-2}
 1-2(1-\delta) \cos(\lambda + \lambda_0) + (1-\delta)^2 \geq  1-c^2_{\max}.
 \end{equation}
From  (\ref{eq.B-1}) and (\ref{eq.B-2}) we obtain
 \begin{equation}
  \frac{f_{U, \text{min}}}{4}  \int^{\pi}_{0} |1-(1-\delta)   e^{-i(\lambda\ -\lambda_0)}|^{-2}d\lambda \leq  \int^{\pi}_{0} f_{\delta}(\lambda)\ d\lambda
   \leq  \frac{f_{U, \text{max}}}{1-c^2_{\max}}    \int^{\pi}_{0} |1-(1-\delta)e^{-i(\lambda\ -\lambda_0)}|^{-2} d\lambda,
\end{equation}
where $ f_{U, \text{min}} =\inf_{\lambda\in[0,\pi] }f_U(\lambda)$ and $ f_{U, \text{max}} =\sup_{\lambda\in[0,\pi] }f_U(\lambda)$. 
Furthermore, 
 \begin{align*}
\int^{\pi}_{0} |1-(1-\delta)   e^{-i(\lambda\ -\lambda_0)}|^{-2} d\lambda
 &= \int^{\pi - \lambda_0 }_{-\lambda_0} |1-(1-\delta)   e^{-i x}|^{-2}dx\\
 & = \int^{\pi - \lambda_0 }_{-\lambda_0} (\delta^2-2(1-\delta)(\cos(\lambda)-1))^{-1} d\lambda\\
 & = \int^{\tan(\frac{\pi - \lambda_0}{2})}_{\tan({-\frac{\lambda_0}{2})}} (\delta^2-2(1-\delta)(\frac{1-t^2}{1+t^2}-1))^{-1}\frac{2}{1+t^2}dt\\
 & = \frac{2}{\delta^2}\int^{\tan(\frac{\pi - \lambda_0}{2})}_{\tan({-\frac{\lambda_0}{2})}} (1 + (\frac{4 (1-\delta)}{\delta^2}+1)t^2)^{-1} dt \\
 & = \frac{2}{\delta(2-\delta)} \Big\{\text{arctan}\Big(\frac{(2-\delta) \tan((\pi - \lambda_0)/2)}{\delta}\Big) \\ & \ \ \ \ \ \ \ \ \  -\text{arctan}\Big(\frac{(2-\delta) \tan(-\lambda_0/2)}{\delta}\Big)\Big\}\\
& = {\mathcal O}(\delta^{-1}),
 \end{align*}
 since $ \lim_{x\rightarrow\pm \infty} {\rm arctan}(x)=\pm\pi/2$. Notice that to obtain the third equality above,  we used the substitution  
$  t=\tan(\lambda/2) $   from which we get $ dt = 1/(2\cos^2(\lambda/2)) d\lambda$, that is,
\[ d\lambda= 2\cos^2(\lambda/2) dt = (1+\cos(\lambda)) dt = \frac{2}{1+t^2} dt,\]
where the last equality  follows using  $ \cos(\lambda) =  (1-t^2)/(1+t^2)  $.
Hence we have shown that 
\begin{equation*}
    \frac{f_{U, \text{min}}}{4} \frac{1}{\delta} \leq \gamma_{\delta}(0) \leq \frac{f_{U, \text{max}}}{1-c^2_{\max}} \frac{1}{\delta}
\end{equation*}
from which we conclude that $ 
     \lim_{\delta\rightarrow 0}\gamma_{\delta}(0) =   + \infty $. 
     
\noindent    Consider (iii). 
Let 
\[ 
g_{peak}(\lambda) = |1-(1-\delta)   e^{-i(\lambda -\lambda_0)}|^{-2}|1-(1-\delta)   e^{-i(\lambda +\lambda_0)}|^{-2}.
\] 
Then, 
 \begin{align*}
 g_{peak}(\lambda) &= \Big|\sum_{k=0}^{\infty} (1-\delta)^k     e^{-ik(\lambda -\lambda_0)}\Big| ^{2} 
 \frac{1}{1-2(1-\delta)  \cos(\lambda +\lambda_0)+(1-\delta)^{2}}\\
  &= \sum_{k_1, k_2 =0}^{\infty} (1-\delta)^{k_1+k_2}     e^{-i({k_1-k_2})(\lambda -\lambda_0)} \frac{1}{1-2(1-\delta)   
  \cos(\lambda +\lambda_0)+(1-\delta)^{2}}\\
  &= \sum_{k_1, k_2 =0}^{\infty} g(k_1, k_2))\ v_{\delta}({\lambda}) ,
  \end{align*}
where $g(k_1, k_2) = (1-\delta)^{k_1+k_2}     e^{-i({k_1-k_2})(\lambda-\lambda_0)}$ and $v_{\delta}({\lambda}) = (1-2(1-\delta)\cos(\lambda+\lambda_0)+(1-\delta)^{2})^{-1}$. We have, 
\begin{align*}
   {\sum_{k_1, k_2 =0}^{\infty}} g(k_1, k_2)\ &= \sum_{k =0}^{\infty}\sum_{m =0}^{\infty} g(m+k, m) + \sum_{k =-\infty}^{-1}\sum_{m =0}^{\infty} g(m, m + |k|)\\
    &= \sum_{k =0}^{\infty}\sum_{m =0}^{\infty} (1-\delta)^{2m+k}     e^{-i{k}(\lambda-\lambda_0)}+ \sum_{k =-\infty}^{-1}
    \sum_{m =0}^{\infty} (1-\delta)^{2m+|k|} e^{-i{k}(\lambda-\lambda_0)} \\
    & = \sum_{m =0}^{\infty} (1-\delta)^{2m} \sum_{k \in \mathbb{Z}} (1-\delta)^{|k|} e^{-i{k}(\lambda-\lambda_0)}.
   \end{align*}
   Hence, 
\begin{equation*}
    g_{peak}(\lambda) = v_{\delta}({\lambda})\sum_{m =0}^{\infty} (1-\delta)^{2m} \sum_{k \in \mathbb{Z}} (1-\delta)^{|k|} e^{-i{k}(\lambda-\lambda_0)}.
\end{equation*}
For  the autocovariance at lag $h$ we have 
\begin{align*}
\gamma_{\delta}(h)\ &= \int^{\pi}_{-\pi} \cos(h\lambda) f_{\delta}(\lambda) d\lambda \\
&= \int^{\pi}_{-\pi} \frac{1}{2} (e^{ih\lambda}+ e^{-ih\lambda}) g_{peak}(\lambda) f_{U}(\lambda) d\lambda \\
&= \frac{1}{2}\int^{\pi}_{-\pi}  (e^{ih\lambda}(\sum_{m =0}^{\infty} (1-\delta)^{2m} \sum_{k \in \mathbb{Z}} (1-\delta)^{|k|} e^{-i{k}(\lambda-\lambda_0)}) v_{\delta}({\lambda}) f_{U}(\lambda) d\lambda \\ 
&+\frac{1}{2}\int^{\pi}_{-\pi}  (e^{-ih\lambda}(\sum_{m =0}^{\infty} (1-\delta)^{2m} \sum_{k \in \mathbb{Z}} (1-\delta)^{|k|} e^{-i{k}(\lambda-\lambda_0)}) v_{\delta}({\lambda}) f_{U}(\lambda) d\lambda\\
& = \frac{1}{2} \sum_{m =0}^{\infty} (1-\delta)^{2m} \sum_{k \in \mathbb{Z}} (1-\delta)^{|k|} e^{i{k}\lambda_0}[\int^{\pi}_{-\pi} e^{-i{(k-h)}\lambda} v_{\delta}({\lambda}) f_{U}(\lambda) d\lambda \\
& \ \ \ \ \ \ \ \ \ \ \ \ \ \ \ \ \ \ \  \ \ \ \ \ \ \ \ \ \ \ \ \ \ \ \ \ \ \ \ \ \ \ \ \ \   + \int^{\pi}_{-\pi} e^{-i{(k+h)}\lambda} v_{\delta}({\lambda}) f_{U}(\lambda) d\lambda]\\
&= \frac{1}{2(1-(1-\delta)^2)} \sum_{k \in \mathbb{Z}} (1-\delta)^{|k|} e^{i{k}\lambda_0} (W(k-h) + W(k+h))
\end{align*}
where $W(s) = \int^{\pi}_{-\pi}  v_{\delta}({\lambda}) f_{U}(\lambda) e^{-i{s}\lambda}d\lambda$.\\

\noindent Note that for the variance we get  
\begin{equation*}
    \gamma_{\delta}(0) = \frac{1}{1-(1-\delta)^2} \sum_{k \in \mathbb{Z}} (1-\delta)^{|k|} e^{i{k}\lambda_0} W(k)
\end{equation*}
Now, we have that  
\begin{align*}
    \sum_{k \in \mathbb{Z}} (1-\delta)^{|k|} e^{i{k}\lambda_0} W(k+h) &= \sum_{l \in \mathbb{Z}} e^{i({l-h})\lambda_0} (1-\delta)^{|l-h|} W(l)\\
    & = e^{-i{h}\lambda_0} \sum_{l \in \mathbb{Z}} e^{il\lambda_0} (1-\delta)^{|l-h|} W(l)
\end{align*}
Similarly, 
\begin{align*}
    \sum_{k \in \mathbb{Z}} (1-\delta)^{|k|} e^{i{k}\lambda_0} W(k-h) =  e^{i{h}\lambda_0} \sum_{l \in \mathbb{Z}} e^{i{l}\lambda_0} (1-\delta)^{|l+h|} W(l)
\end{align*}
Also note that $||a|-|b||\leq |a-b|$  which implies $ -|a-b|\leq |a|-|b|\leq |a-b|$ and that  if $ v\leq w$ then $(1-\delta)^w \leq (1-\delta)^v$. we then obtain,   
\begin{align*}
&|\sum_{l \in \mathbb{Z}} e^{il \lambda_0} (1-\delta)^{|l|} W(l)-\sum_{l \in \mathbb{Z}} e^{il\lambda_0} (1-\delta)^{|l \pm h|} W(l)| \\ &\leq \sum_{l \in \mathbb{Z}} |(1-\delta)^{|l \pm h|} - (1-\delta)^{|l|}| \,
|W(l)| \\
&\leq \sum_{l \in \mathbb{Z}} (1-\delta)^{|l|} \text{max} \{|1-(1-\delta)^{|h|}|,|(1-\delta)^{-|h|}-1|\} |W(l)| =\mathcal{O}(\delta)
\end{align*}
Note that this is valid because $W(l)$ is absolutely summable, since the $W(l)$ represent the Fourier coefficients of the differentiable function $v_{\delta}({\lambda}) f_{U}(\lambda), \lambda \in [-\pi, \pi]$.\\
Hence
\begin{align*}
\frac{\sum_{l \in \mathbb{Z}} e^{il\lambda_0} (1-\delta)^{|l \pm h|} W(l)} {\sum_{l \in \mathbb{Z}} e^{il \lambda_0} (1-\delta)^{|l|} W(l)} - 1  & =  \frac{\sum_{l \in \mathbb{Z}} e^{il\lambda_0} (1-\delta)^{|l \pm h|} W(l) - \sum_{l \in \mathbb{Z}} e^{il \lambda_0} (1-\delta)^{|l|} W(l)} {\sum_{l \in \mathbb{Z}} e^{il \lambda_0} (1-\delta)^{|l|} W(l)}  =\mathcal{O}(\delta)
\end{align*}
since
\begin{align*}
    \min_{\delta \in (0,1)} \sum_{l \in \mathbb{Z}} e^{il \lambda_0} (1-\delta)^{|l|} W(l)  & \ge   \min_{\delta \in (0,1)}\{ (1-(1-\delta)^2) \gamma_{\delta}(0)\}
    & \ge \min_{\delta \in (0,1)}\{  \delta(2 - \delta) \frac{f_{U, min}}{4} \frac{1}{\delta}\} 
    & \ge \frac{f_{U, \min}}{4}.
\end{align*}
Then,
\begin{align*}
\frac{\sum_{l \in \mathbb{Z}} e^{il\lambda_0} (1-\delta)^{|l \pm h|} W(l)} {\sum_{l \in \mathbb{Z}} e^{il \lambda_0} (1-\delta)^{|l|} W(l)} = 1 + \mathcal{O}(\delta).
\end{align*}

\noindent Therefore  for $\ h\in \mathbb{N}_0$, 
\begin{align*}
\rho_{\delta}(h)\ &= \frac{\gamma_{\delta}(h)}{\gamma_{\delta}(0)}\\
& = \frac{\frac{1}{2(1-(1-\delta)^2)} [e^{+ih\lambda_0}\sum_{l \in \mathbb{Z}} (1-\delta)^{|l+h|} e^{i{l}\lambda_0} W(l)+ e^{-ih\lambda_0}\sum_{l \in \mathbb{Z}} (1-\delta)^{|l-h|} e^{i{l}\lambda_0} W(l)]}
{\frac{1}{(1-(1-\delta)^2)} \sum_{k \in \mathbb{Z}} (1-\delta)^{|k|} e^{i{k}\lambda_0} W(k)} \\
&=   \frac{1}{2} (e^{+ih\lambda_0} +e^{-ih\lambda_0})  + {\mathcal O}(\delta) = \cos(h\lambda_0) + {\mathcal O}(\delta), \ \ \mbox{as $ \delta \rightarrow 0$}.
\end{align*}
To see (iv) note that $A_\delta(z)=1-2(1-\delta)\cos(\lambda_0)z +(1-\delta)^2z^2$ has the roots $z_{1,2}= (1-\delta)^{-1}e^{\pm i\lambda_0}$. From this  we   get  after some simple  algebra  that  ($L$ denotes the usual lag operator)
\begin{align*}
\big(1-2(1-\delta)\cos(\lambda_0)L & +(1-\delta)^2L^2\big)^{-1}  = 1+2(1-\delta) \cos(\lambda_0) L \\
&  + \sum_{n=2}^\infty (1-\delta)^n \big[2\cos(n\lambda_0) +\sum_{k=1}^{n-1}\cos((2k-n)\lambda_0)\big] L^n.
\end{align*}
which implies  that  for $ j>2$ and even,   $ b_{j,\delta} =(1-\delta)^j\big(1+  2\cos(j\lambda_0) + 2\sum_{s=1}^K \cos(2s\lambda_0)\big)$,  where $ K=\lfloor j/2\rfloor -1$, while 
for $ j>2$ and odd,    $ b_{j,\delta} =(1-\delta)^j\big(   2\cos(j\lambda_0) + 2\sum_{s=1}^{K+1} \cos((2s-1)\lambda_0)\big)$.
The assertion then follows by straightforward calculations  and using the relations $ \cos(2x)=2\cos^2(x) -1 $ and $\sum_{s=1}^n \sin^2(x)=n/2 -\cos((n+1)x)\sin(nx)/2\sin(x)$ for  $ \sin(x)\neq 0$. \hfill $\Box$\\

\noindent {\bf Proof of Lemma~\ref{le.2}:}
(i) Fitting an AR(2)-model by means of  Yule-Walker estimators leads to the coefficients 
\[ 
(a_{1,\delta}, a_{2,\delta}) = \Big(\frac{1}{1-\rho_{1,\delta}^2}(\rho_{1,\delta}(1-\rho_{2,\delta})), \  \frac{1}{1-\rho_{1,\delta}^2} (\rho_{2,\delta}-\rho^2_{1,\delta}) \Big) 
\]
from which assertion (i) follows using Lemma 1(iii). \\
Since for   $ a_{2,\delta}<0$,  
$ \lambda_{0,\delta} = {\rm arc cos}(a_{1,\delta}(1-a_{2,\delta})/(-4a_{2,\delta}))$,  
we immediately get using (i) and  for $ \delta\rightarrow0$, that  
$ \lambda_{0,\delta} \rightarrow  {\rm arc cos}(\cos(\lambda_0))=\lambda_0$, which is assertion (ii). \hfill $\Box$\\

\noindent {\bf Proof of Proposition~\ref{prop.1}:} \  Let $ {\mathcal F}(n) =\{ -\lfloor (n-1)/2 \rfloor, \ldots, [n/2]\}$ be the set of Fourier frequencies. For $h\in\{0,1,\ldots,n-1\}$ and $ \lambda_{j,n} =2\pi j/n$ we have  using 
$  I_{n,\delta}(\lambda_{j,n}) = (2\pi)^{-1}\sum_{|s|<n} \widetilde{\gamma}_\delta(s) e^{-i s \lambda_{j,n}}$,  that 
\begin{align*}
\frac{2\pi}{n} \sum_{j \in {\mathcal F}(n)} I_{n,\delta}(\lambda_{j,n}) e^{i\lambda_{j,n} h} & =  \frac{1}{n}\sum_{|s|<n} \widetilde{\gamma}_\delta(s) \sum_{j\in{\mathcal F}(n)}\Big(e^{i 2\pi (h-s)/n}\Big)^j.
\end{align*}
For  $ s=h$ or $ s=-(n-h)$  we have  $\sum_{j\in{\mathcal F}(n)}\Big(e^{i 2\pi (h-s)/n}\Big)^j =n$,  since  $ |{\mathcal F}(n)|=n$.  For  $ s \notin \{h, -(n-h)\}$ let $\varrho:=e^{i2\pi(h-s)/n} $ and note that $\varrho \neq 1$. Therefore,
\begin{align*}
\sum_{j\in{\mathcal F}(n)}\Big(e^{i 2\pi (h-s)/n}\Big)^j & = \frac{\displaystyle  \varrho^{\lfloor n/2\rfloor +1} - \varrho^{-\lfloor (n-1)/2\rfloor}}{\displaystyle \varrho -1}\\
& = \left\{ \begin{array}{lll} \frac{\displaystyle \varrho}{\displaystyle \varrho -1} \big( \varrho^{n/2} - \varrho^{-n/2}\big)   & & \mbox{for  $n$ even} \\
& & \\
 \frac{\displaystyle \varrho^{1/2}}{\displaystyle \varrho -1} \big( \varrho^{n/2} - \varrho^{-n/2}\big)   & & \mbox{for  $n$ odd} 
 \end{array}\right.
\end{align*} 
and   in both cases,  $ \sum_{j\in{\mathcal F}(n)}\Big(e^{i 2\pi (h-s)/n}\Big)^j=0$ since $ \varrho^{n/2} - \varrho^{-n/2} =-2i \sin(\pi s)=0$. This implies for  $\widehat{\gamma}_\delta(h)$ and $h=0$, that 
\begin{align*}
\widehat{\gamma}_\delta(0)&  =\widetilde{\gamma}_\delta(0) + \frac{2\pi}{n} \Big( \sum_{j\in{\mathcal G}(n)} I_{n,\delta}(\lambda_{j,n}) -\sum_{j\in{\mathcal F}(n) }I_{n,\delta}(\lambda_{j,n})  \Big)\\
 &  =\widetilde{\gamma}_\delta(0) - n^{-1}\big( I_{n,\delta}(0) - \mathds{1}(n \ \mbox{even}) I_{n,\delta}(\pi) \big)\\
 & =\widetilde{\gamma}_\delta(0)  + O_P(n^{-1}\delta^{-2}),
\end{align*}
since  for all $ \lambda \in[0,\pi]$,
\begin{align} \label{eq.E-Per}
{\rm E}(I_{n,\delta}(\lambda) ) 
& = \frac{1}{2\pi n} \sum_{j_1=0}^\infty\sum_{j_2=0}^\infty b_{j_1,\delta} b_{j_2,\delta}  
\sum_{|s|,|t|<n}{\rm E}(U_{t-j_1}U_{s-j_2})e^{-i\lambda(t-s)} \nonumber \\
& \leq \frac{1}{2\pi } \Big( \sum_{j=0}^\infty |b_{j,\delta}|\Big)^2 \sum_{h\in \mathbb{Z}}|\gamma_U (h)| 
= O(\delta^{-2}), 
\end{align}
since  $ |b_{j,\delta}|   \leq C (1-\delta)^j $; see Lemma~\ref{le.1}(iv). 

For $ h \in \{1,2, \ldots, n-1\}$ we get 
\begin{align*}
\widehat{\gamma}_\delta(h)&  =\widetilde{\gamma}_\delta(h) + \widetilde{\gamma}_\delta (n-h)  + \frac{2\pi}{n} \Big( \sum_{j\in{\mathcal G}(n)} I_{n,\delta}(\lambda_{j,n}) e^{-i\lambda_{j,n} h}
-\sum_{j\in{\mathcal F}(n) }I_{n,\delta}(\lambda_{j,n})e^{-i\lambda_{j,n} h } \Big)\\
&  =\widetilde{\gamma}_\delta(0) +  \widetilde{\gamma}_\delta (n-h) - n^{-1}\big( I_{n,\delta}(0) -  
\mathds{1}(n \ \mbox{even}) I_{n,\delta}(\pi)e^{i\pi h}\big) \\
 &  =\widetilde{\gamma}_\delta(h) + O_P(n^{-1}\delta^{-2}),
\end{align*}
where the last equality follows because for any fixed $h \in \{1,2, \ldots, n-1\}$,
\begin{align*}
{\rm E}|\widetilde{\gamma}_\delta (n-h)| & \leq \frac{1}{n} \sum_{t=1}^{h} {\rm E}|X_{t+h,\delta}X_{t,\delta}| \\
& \leq \frac{1}{n} \sum_{j_1=0}^\infty \sum_{j_2=0}^\infty |b_{j_1,\delta}||b_{j_2,\delta}| \,  
    \sum_{t=1}^{h} {\rm E}|U_{t+h-j_1} U_{t-j_2} | \\
& \leq  \frac{h}{n} \Big(\sum_{j=0}^\infty b_{j,\delta}\Big)^2 \gamma_U(0) = O(n^{-1}\delta^{-2}),
\end{align*}
and   $ I_{n,\delta}(0) +  \mathds{1}(n \ \mbox{even}) I_{n,\delta}(\pi)e^{i\pi h} = O_P(\delta^{-2})$; see  (\ref{eq.E-Per}).
\hfill $\Box$\\

\noindent {\bf Proof of Lemma~\ref{le.Sum}:}
We proof (ii) since (i) is a direct consequence of (ii).  Using  the expression 
      $ \gamma_{\delta}(h) =  \sum_{s = -\infty}^{\infty}\sum_{j = 0}^{\infty} b_{j,\delta} b_{j + |s|,\delta} \gamma_U(h-s)$,
we get
     \begin{align*}
        \sum_{h \in \mathbb{Z}}|h||\gamma_{\delta}(h)| 
        &\leq    \sum_{h \in \mathbb{Z}}\sum_{s = -\infty}^{\infty}\sum_{j = 0}^{\infty} |b_{j,\delta}| |b_{j+ |s|,\delta}| |h| |\gamma_U(h-s)|\\
        & =  \sum_{s = -\infty}^{\infty}\sum_{j = 0}^{\infty} |b_{j,\delta}| |b_{j + |s|,\delta}| \sum_{h \in \mathbb{Z}}|h| |\gamma_U(h-s)|\\
        & \leq C\big(2\sum_{s=0}^\infty (1-\delta)^s-1\big) \sum_{j=0}^\infty (1-\delta)^{2j}\\
        & = \mathcal{O} (\delta^{-2}).
        \end{align*}
To see (iii) recall that 
\begin{align*}
\textrm{cum}(X_{0,{\delta_n}}, & X_{h_1,{\delta_n}}, X_{h_2,{\delta_n}}, X_{h_3,{\delta_n}} ) \\
&  = \sum_{j_0 = 0}^{\infty} \sum_{j_1 = 0}^{\infty} \sum_{j_2 = 0}^{\infty} \sum_{j_3 = 0}^{\infty} b_{j_0,\delta}b_{j_1,\delta}b_{j_2,\delta}b_{j_3,\delta}\, \textrm{cum}(U_{-{j_0}}, U_{{h_1}-{j_1}}, U_{{h_2}-{j_2}},U_{{h_3}-{j_3}} )\\
& = \sum_{j_0 = 0}^{\infty} \sum_{j_1 = 0}^{\infty} \sum_{j_2 = 0}^{\infty} \sum_{j_3 = 0}^{\infty} b_{j_0,\delta}b_{j_1,\delta}b_{j_2,\delta}b_{j_3,\delta}\, \textrm{cum}(U_{0}, U_{{h_1}+ ({j_0}-{j_1})}, U_{{h_2}+ ({j_0}-{j_2})},U_{{h_3}+ ({j_0}-{j_3})} )
\end{align*}
Therefore,
\begin{align*}
& \sum_{h_1, h_2, h_3 = -\infty}^{\infty} (1+|h_1|+|h_2|+|h_3|) |\textrm{cum}(X_{0,{\delta_n}}, X_{h_1,{\delta_n}}, X_{h_2,{\delta_n}}, X_{h_3,{\delta_n}} )|  \\
& \le \sum_{h_1, h_2, h_3 = -\infty}^{\infty} \sum_{j_0, j_1, j_2, j_3= 0}^{\infty} (1+|h_1|+|h_2|+|h_3|) |b_{j_0,\delta}||b_{j_1,\delta}||b_{j_2,\delta}||b_{j_3,\delta}|\\
&\ \ \ \ \ \ \ \ \ \ \  \  \times |\textrm{cum}(U_{0}, U_{{h_1}+ ({j_0}-{j_1})}, U_{{h_2}+ ({j_0}-{j_2})},U_{{h_3}+ ({j_0}-{j_3})} )|\\
& = \sum_{j_0, j_1, j_2, j_3= 0}^{\infty} |b_{j_0,\delta}||b_{j_1,\delta }||b_{j_2,\delta }||b_{j_3,\delta }| \\
& \ \ \ \ \ \ \times \sum_{h_1, h_2, h_3 = -\infty}^{\infty} (1+|h_1|+|h_2|+|h_3|) |\textrm{cum}(U_{0}, U_{{h_1}+ ({j_0}-{j_1})}, U_{{h_2}+ ({j_0}-{j_2})},U_{{h_3}+ ({j_0}-{j_3})} )|\\
&\le C \Big\{ \sum_{j= 0}^{\infty} (1-\delta)^j\Big\}^4=  \mathcal{O}(\delta^{-4}) .   \\
& \hspace{14cm}  \Box 
\end{align*}


\noindent {\bf Proof of Lemma~\ref{le.Per}:}\ First, consider (i). We have, 
\begin{align*}
    {\rm E}[I_{n,\delta}(\lambda_j)] 
    &= \frac{1}{2\pi} \sum_{h=-n+1}^{n-1} (1-\frac{|h|}{n})\ \gamma_{\delta}(h)\ e^{-i \lambda_j h}
\end{align*}
and therefore, 
\begin{align*}
    \big|{\rm E}[I_{n,\delta}(\lambda_j)] - f_{\delta}(\lambda_j)\big|  
    &\leq \big|\frac{1}{2\pi} \sum_{|h| \ge n}\gamma_{\delta}h\ e^{-i \lambda_j(h)}\big| + \big|\frac{1}{2\pi} \sum_{h=-n+1}^{n-1} \frac{|h|}{n} \gamma_{\delta}(h)\ e^{-i \lambda_j h}\big|\\
\end{align*}
From   Lemma~\ref{le.Sum} we get 
%
   $ \big|(2\pi n)^{-1} \sum_{h=-n+1}^{n-1} |h| \gamma_{\delta}(h)\ e^{-i \lambda_j h}\big|= \mathcal{O}(n^{-1}\delta^{-2}) $ and
\begin{align*}
   \Big|\frac{1}{2\pi} \sum_{|h| \ge n}\gamma_{\delta}(h)\ e^{-i \lambda_jh}\Big|   
   &\leq \Big|\frac{1}{2\pi n} \sum^\infty_{h = -\infty} |h| \big|\gamma_{\delta_n}(h)\Big|
    = \mathcal{O}\Big(\frac{1}{n \delta^2}\Big),
\end{align*}
which implies (i).

\noindent To establish (ii) we borrow arguments from the proof of Lemma 2.1 in  \citeasnoun{meyer_pap2023}.
Notice first that,
\begin{align*}
\cov [I_{\delta_n,n}(\lambda_{k,n}), I_{n,\delta}(\lambda_{m,n})] & = \frac{1}{4\pi^2 n^2}\sum_{t_1, t_2, t_3, t_4 = 1}^n \gamma_{\delta}(t_1-t_3)\gamma_{\delta}(t_2-t_4)e^{-i(t_1-t_2)\lambda_{k,n} +  i(t_3-t_4) \lambda_{m,n}}\\
&+\frac{1}{4\pi^2 n^2}\sum_{t_1, t_2, t_3, t_4 = 1}^n \gamma_{\delta}(t_1-t_4)\gamma_{\delta}(t_2-t_3)e^{-i(t_1-t_2)\lambda_{k,n} +  i(t_3-t_4) \lambda_{m,n}}\\
&+\frac{1}{4\pi^2 n^2}\sum_{t_1, t_2, t_3, t_4 = 1}^n \textrm{cum}(X_{t_1,\delta}, X_{t_2,\delta}, X_{t_3,\delta}, X_{t_4,\delta}) \\
& \ \ \ \ \ \ \times e^{-i(t_1-t_2)\lambda_{k,n} +  i(t_3-t_4) \lambda_{m,n}}\\
&= S_1+S_2+S_3,
\end{align*}
with an obvious notation for $ S_i$,  $i=1,2,3$. We consider each of the terms separately.\\

\noindent We begin with the term $S_1$. To handle  this term, we define the function  
$ g_n (\lambda) = \sum_{k=1}^n e^{ik\lambda}$. For  $ \lambda \notin 2\pi \Z$ we have by the geometric sum formula, 
\begin{align*}
     g_n (\lambda) &= e^{i\lambda} \frac{e^{in\lambda}-1}{e^{i\lambda}-1}
     = e^{i(n+1)\frac{\lambda}{2}} \frac{\sin(\frac{n\lambda}{2})}{\sin(\frac{\lambda}{2})}.
\end{align*}
For $\lambda \in \R$, let 
$ g_n (\lambda) g_n (-\lambda) = 2 \pi n F_n(\lambda)$ with $ F_n$ the Fejer kernel defined as 
\begin{align*}
    F_n(\lambda) =  \left\{ \begin{array}{lll} 
\frac{ 1}{ 2 \pi n}\sin^2(n\lambda/2) \big/\sin^2(\lambda/2) & & \mbox{if} \ \  \lambda  \notin 2\pi \mathbb{Z}\\
& & \\
n/2 \pi \ & & \mbox{if}\ \   \lambda  \in 2\pi \mathbb{Z}.
\end{array}  \right.
\end{align*}
 Let 
 $$ z_{n, \lambda_1, \lambda_2} (\lambda) = g_n (\lambda-\lambda_1) g_n (-(\lambda-\lambda_2)). $$
 The $\ell$-th Fourier coefficient  $\widehat{z}_{n, \lambda_1, \lambda_2} [\ell]$ of $ z_{n, \lambda_1, \lambda_2} (\cdot)$  is given by 
 \begin{align*}
 \widehat{z}_{n, \lambda_1, \lambda_2} [\ell] &= \frac{1}{2\pi} \int^\pi_{-\pi} g_n (\lambda-\lambda_1) g_n (-(\lambda-\lambda_2)) e^{-i \ell \lambda} d\lambda\\
 &= \frac{1}{2\pi} \sum^n_{p,q=1} e^{-ip\lambda_1} e^{iq\lambda_2} \int^\pi_{-\pi}  e^{-i(p-q-\ell)\lambda} d\lambda\\
 &= \sum^n_{p,q=1} e^{-ip\lambda_1} e^{iq\lambda_2} \mathds{1}_{\{p-q=\ell\}}.
 \end{align*}
Observe that   $\hat{z}_{n, \lambda_1, \lambda_2} [h] = 0\ \forall\ |h| \ge n$.
Using, 
 $ \gamma_{\delta}(h) \ = \int^{\pi}_{-\pi} e^{ih \lambda}  f_{\delta}(\lambda) d\lambda $, we can write  
 \begin{align} \label{eq.S1-0}
   S_1 & =  \frac{1}{4\pi^2 n^2} \int^{\pi}_{-\pi} z_{n, \lambda_{k,n}, \lambda_{m,n}} (\lambda) f_{\delta}(\lambda) d\lambda  \int^{\pi}_{-\pi} z_{n, -\lambda_{k,n}, -\lambda_{m,n}} (\lambda) f_{\delta}(\lambda) d\lambda. 
 \end{align}  
 We treat the cases $  \lambda_{k,n}=\lambda_{m,n}$ and $ \lambda_{k,n}\neq \lambda_{m,n} $ separately. 
 
\noindent  Case $ \lambda_{k,n}=\lambda_{m,n}$:  For this case  we have,  
 \begin{equation} \label{eq.S1-1}
  S_1=  \int^{\pi}_{-\pi} F_{n} (\lambda-\lambda_{k,n})f_{\delta}(\lambda) d\lambda  \int^{\pi}_{-\pi} F_{n} (\lambda+\lambda_{k,n}) f_{\delta}(\lambda) d\lambda.
  \end{equation}
 It yields $ F_n(\cdot -\lambda_{k.n}) =z_{n,\lambda_{k,n},\lambda_{k,n}}(\cdot)/(2\pi n)$ and the $\ell$-th Fourier coefficient 
 of $ F_{n} (\cdot -\lambda_{k,n})$ is given by 
\begin{align*}
\widehat{F}_{n} (\cdot -\lambda_{k,n})[\ell] &= \widehat{z}_{n, \lambda_{k,n}, \lambda_{k,n}} [\ell] = \Big(\frac{n-|\ell|}{2\pi n}\Big)\ \ e^{-i \lambda_{k,n} \ell}\mathbf {1}_{\{|\ell|< n\}}.
\end{align*}
Furthermore,  the $\ell$-th  Fourier coefficient of $f_{\delta}(\lambda)$ is given by
\begin{align*}
\widehat{f}_{\delta}[\ell]  &= \frac{1}{2\pi} \int^{\pi}_{-\pi} e^{-i \ell  \lambda}  f_{\delta}(\lambda) d\lambda =  \frac{\gamma_{\delta}(-\ell)}{2\pi}
\end{align*}
By Parseval's Theorem, we get that
\begin{align*}
\int^{\pi}_{-\pi} F_{n} (\lambda-\lambda_{k,n})f_{\delta}(\lambda) d\lambda &= 2\pi \sum_{\ell=-\infty}^{\infty} \widehat{f}_{\delta}(\ell) \overline{\widehat{F}_{n} (\cdot -\lambda_{k,n})[\ell]}\\
&= \frac{1}{2\pi} \sum_{\ell=-n+1}^{n-1} \Big(1-\frac{|\ell|}{n}\Big)\ \gamma_{\delta}(\ell)\ e^{-i \lambda_k \ell}
\end{align*}
and as in the proof of assertion (i),
\[ \big|\int^{\pi}_{-\pi} F_{n} (\lambda-\lambda_{k,n})f_{\delta}(\lambda) d\lambda - f_{\delta}(\lambda_{k,n})  \big|= {\mathcal O}\big(n^{-1}\delta^{-2}\big) \]
Using analogous arguments, we obtain that   
\[ 
\Big| \int^{\pi}_{-\pi} F_{n} (\lambda+\lambda_{m,n}) f_{\delta}(\lambda) d\lambda -f_{\delta}(\lambda_{k,n}) \Big| 
= {\mathcal O}\big(n^{-1}\delta^{-2}\big) ,
\] 
from which we  get in view of  (\ref{eq.S1-1}) and $ f_{\delta}(\lambda) = |A_{\delta}(e^{-i\lambda})|^{-2}f_U(\lambda)$, that   
\begin{align*}
S_1 &=   f^2_{\delta}(\lambda_{k,n}) + \mathcal{O}\Big(\frac{1}{|A_{\delta}(e^{-i\lambda_{k,n}})|^2 n\delta^2}\Big) +  \mathcal{O}\Big(\frac{1}{n^2\delta^4}\Big).
\end{align*}

\noindent Case $\lambda_{k,n} \neq \lambda_{m,n}$:   
In this case,  $\lambda_{k,n}-\lambda_{m,n}=2\pi(k-m)/n$ is a Fourier coefficient  different from $ 2\pi \Z$ and therefore, $\sum_{q=1}^n e^{i q(\lambda_{k,n}-\lambda_{m,n}) }=0 $. Using this we obtain the bound 
for the Fourier coefficients of $z_{n,\lambda_{k,n},\lambda_{m,n}}$ for $|\ell|<n$, 
\begin{align*}
    |\widehat{z}_{n, \lambda_{k,n}, \lambda_{m,n}} [\ell] | &=  |\sum^{n-|\ell|}_{q=1} e^{-iq(\lambda_{m,n} - \lambda_{k,n})}|\\
    &= |-\sum^{n}_{q=n-|\ell|+1} e^{-iq(\lambda_{m,n} - \lambda_{k,n})}| \leq |\ell|.
\end{align*}
Recall that the Fourier coefficients of $z_{n,\lambda_{k,.n},\lambda_{m,n}}$  for $ \ell \geq n$ are zero. 
From this and Parseval's Theorem we get 
\begin{align*}
\Big|\int^{\pi}_{-\pi} z_{n, \lambda_{k,n}, \lambda_{m,n}} (\lambda)  f_{\delta}(\lambda) d\lambda \Big| 
&= \Big| 2\pi\sum_{\ell=-\infty}^{\infty} \widehat{z}_{n, \lambda_{k,n}, \lambda_{m,n}} [\ell] 
\overline{\widehat{f}_{\delta}(\ell)}\Big| \\
& \leq \sum_{\ell=-\infty}^{\infty} (1+|\ell|) |\gamma_{\delta}(\ell)| = \mathcal{O} \Big(\frac{1}{\delta^2}\Big),
\end{align*}
by Lemma~\ref{le.Sum}.
Hence,  in this case and in view of (\ref{eq.S1-0}) we have $ S_1 ={\mathcal  O}(n^{-2}\delta^{-4})$.\\

\noindent Consider next  the term $ S_2$ for which we get,
$$ S_2 =  \frac{1}{4\pi^2 n^2} \int^{\pi}_{-\pi} z_{n, \lambda_{k,n}, -\lambda_{m,n}} (\lambda) f_{\delta_n}(\lambda) d\lambda  \int^{\pi}_{-\pi} z_{n, -\lambda_{k,n}, \lambda_{m,n}} (\lambda) f_{\delta_n}(\lambda) d\lambda.$$
From this we have  if $ \lambda_{k,n}=\lambda_{m,n}=0$ that 
\begin{align*}
 S_2 & =\int_{-\pi}^\pi F_n(\lambda)f_{\delta}(\lambda) d\lambda \int_{-\pi}^\pi F_n(\lambda)f_{\delta}(\lambda) d\lambda \\
 & =  f^2_{\delta}(0) + {\mathcal O}\Big(\frac{1}{|A_{\delta}(1)|^2 n \delta^2}\Big) + {\mathcal O}\Big(\frac{1}{n\delta^2} \Big). 
 \end{align*}
 If $ \lambda_{k,n}=\lambda_{m,n}\in \{ 0,\pi\}$ we get arguing as in the proof of (1.6)  in the 
 supplementary file of  \citeasnoun{meyer_pap2023}, that 
 \begin{align*}
 S_2 
 & =  f^2_{\delta}(\pi) + {\mathcal O}\Big(\frac{1}{|A_{\delta}(1)|^2 n \delta^2}\Big) + {\mathcal O}\Big(\frac{1}{n\delta^2} \Big). 
 \end{align*}
For all other cases  of $ \lambda_{k,n}$ and $\lambda_{m,n}$, that is for $ \lambda_{k,n}=\lambda_{m,n} \notin \{0,\pi\}$ or $ \lambda_{k,n}\neq \lambda_{m,n}$,  we get 
\[ \big| \widehat{z}_{n, \lambda_{k,n}, -\lambda_{m,n}} [\ell] \big| \leq \big|\ell\big| \ \ \mbox{and} \ \   \big| \widehat{z}_{n,- \lambda_{k,n}, \lambda_{m,n}} [\ell] \big| \leq \big|\ell\big| ,\]
and arguing as for $ S_1$ we obtain $ S_2={\mathcal O}(n^{-2}\delta^{-4})$.

\noindent  Consider  next  the term $S_3$.  Note first that 
\begin{align*}
   f_{4,X_{\delta_n}} ({\lambda_{k,n}}, {\lambda_{m,n}}, -{\lambda_{m,n}})  
   & = |A_{\delta}(e^{-i{\lambda_{k,n}}})|^{-2} |A_{\delta}(e^{-i{\lambda_{m,n}}})|^{-2}f_{4,U} ({\lambda_{k,n}}, {\lambda_{m,n}}, -{\lambda_{m,n}}).
\end{align*}
Then,
\begin{align*}
\frac{n}{2\pi} S_3 &= \frac{1}{(2\pi)^3 n}\sum_{t_1, t_2, t_3, t_4 = 1}^n \textrm{cum}\big(X_{t_1,\delta}, X_{t_2,\delta}, X_{t_3,\delta}, X_{t_4,\delta} \big)e^{-i(t_1-t_2)\lambda_k +  i(t_3-t_4) \lambda_m}\\ 
&  = \frac{1}{(2\pi)^3 n}\sum_{t_4 = 1}^n\sum_{t_1, t_2, t_3 =1 -t_4}^{n-t_4} \textrm{cum}(X_{t_1-t_4,\delta}, X_{t_2-t_4,\delta}, X_{t_3-t_4,\delta}, X_{0,\delta} )e^{-i(t_1-t_2)\lambda_{k,n} +  i t_3 \lambda_{m,n}}\\
&=  \frac{1}{(2\pi)^3 n}\sum_{h_1, h_2, h_3 = -(n-1)}^{n-1} q_n(h_1, h_2, h_3) \textrm{cum}(X_{h_1,\delta}, X_{h_2,\delta}, X_{h_3,\delta}, X_{0,\delta} )e^{-i (h_1 -h_2)\lambda_{k,n} -  i h_3 \lambda_{m,n}}
\end{align*}
where $q_n(h_1, h_2, h_3)$ counts how often the respective summand  appears. It holds that 
 \[
 q_n(h_1, h_2, h_3) = \big(n-\max\{|h_1|, |h_2|, |h_3|, |h_1 - h_2|, |h_1 - h_3|, |h_2 - h_3|\}\big)_{+} .
 \]
 As in the proof of (1.1) in the supplementary file of  \citeasnoun{meyer_pap2023},
we can replace $q_n(\cdot)$ by $n$, 
 whereas the remainder is bounded by 
 \begin{align*}
 \Big|\frac{1}{(2\pi)^3n} & \sum_{h_1,h_2,h_3=-(n-1)}^{n-1} \big( q_n(h_1,h_2,h_3) - n\big) \textrm{cum}\big(X_{h_1,\delta}, X_{h_2,\delta}, X_{h_3,\delta}, X_{0,\delta} \big) 
 e^{-i (h_1 -h_2)\lambda_{k,n} -  i h_3 \lambda_{m,n}} \Big| \\
 & \leq 
 \frac{1}{(2\pi)^3 n}
 \sum_{h_1,h_2,h_3=-\infty}^{\infty} 2(1+|h_1|+|h_2|+ |h_3|)
 \big|\textrm{cum}\big(X_{h_1,\delta}, X_{h_2,\delta}, X_{h_3,\delta}, X_{0,\delta} \big)\big|\\
 & = {\mathcal O}(n^{-1}\delta^{-4})
 \end{align*}
 by Lemma~\ref{le.Sum}. Hence,
 \begin{align*}
\frac{n}{2\pi} S_3 &=   \frac{1}{(2\pi)^3 }\sum_{h_1, h_2, h_3 = -(\infty)}^{\infty}  \textrm{cum}(X_{h_1,{\delta}}, X_{h_2,{\delta}}, X_{h_3,{\delta}}, X_{0,{\delta}} )e^{-i (h_1 -h_2)\lambda_{k,n} -  i h_3 \lambda_{m,n}}\\
& \ \ \ \  + D_n + {\mathcal O}(n^{-1}\delta^{-4}),
\end{align*}
where  
\begin{align*}
\big|D_n\big| & =\Big| \frac{1}{(2\pi)^3 }\sum_{h_1, h_2, h_3 = -(n-1)}^{n-1}  \textrm{cum}(X_{h_1,{\delta}}, X_{h_2,{\delta}}, X_{h_3,{\delta}}, X_{0,{\delta}} )e^{-i (h_1 -h_2)\lambda_{k,n} -  i h_3 \lambda_{m,n}}\\
& \ \ \ \ \ \  - \frac{1}{(2\pi)^3 }\sum_{h_1, h_2, h_3 = -\infty}^{\infty}  \textrm{cum}(X_{h_1,{\delta}}, X_{h_2,{\delta}}, X_{h_3,{\delta}}, X_{0,{\delta}} )e^{-i (h_1 -h_2)\lambda_{k,n} -  i h_3 \lambda_{m,n}}\Big|\\
& \leq \frac{1}{(2\pi)^3 }\sum_{h_1, h_2, h_3 = -\infty}^{\infty}  {\bf 1}_{\max\{|h_1|, |h_2|,|h_3|\geq n\}}\big|\textrm{cum}(X_{h_1,{\delta}}, X_{h_2,{\delta}}, X_{h_3,{\delta}}, X_{0,{\delta}} )\big|\\
& \leq   \frac{1}{(2\pi)^3 } \frac{1}{n}\sum_{h_1, h_2, h_3 = -\infty}^{\infty} \big(|h_1|+ |h_2|+|h_3|\big)\big|\textrm{cum}(X_{h_1,{\delta}}, X_{h_2,{\delta}}, X_{h_3,{\delta}}, X_{0,{\delta}} )\big|\\
& = {\mathcal O}(n^{-1}\delta^{-4}),
\end{align*}
again by Lemma~\ref{le.Sum}.
Therefore, 
\[  S_3   =\frac{\displaystyle 2\pi}{\displaystyle n|A_{\delta}(\lambda_{k,n})|^2 |A_{\delta}(\lambda_{m,n})|^{2}}  f_{4,U}(\lambda_{k,n}, -\lambda_{m,n}, \lambda_{m,n})    + \mathcal{O}(\frac{1}{n^2\delta^4}).
\hspace*{2cm} \Box
\]
\vspace*{0.5cm}

\noindent {\bf Proof of Lemma~\ref{le.gamma0}:}  Using Lemma~\ref{le.1} and ~\ref{le.Per} we obtain, 
 \begin{align*}
    {\rm E}\frac{\hat{\gamma}_{\delta_n}({0})}{\gamma_{\delta_n}({0})} &= \frac{1}{\gamma_{\delta_n}({0})}\frac{2\pi}{n}\sum_{\lambda_j \in {\mathcal G}_n}{\rm E}\, I_{n,\delta_n} (\lambda_j)\\
     &= \frac{1}{\gamma_{\delta_n}({0})}\frac{2\pi}{n}\sum_{\lambda_j \in {\mathcal G}_n} 
     \big( f_{\delta_n} (\lambda_j) + \mathcal{O}(\frac{1}{n\delta_n^2})\big) \\
      &= \frac{1}{\gamma_{\delta_n}({0})}\Big(
      \int_{-\pi}^{\pi} f_{\delta_n} (\lambda)d\lambda + \mathcal{O}(\frac{1}{n\delta_n^2})\Big)  \\
      &=  1+\mathcal{O}((n\delta_n)^{-1})  \xrightarrow[\delta_n \to 0]{} 1.
\end{align*}
Furthermore,  as in the proof of Lemma~\ref{le.1}(ii),  it holds true that 
\begin{equation} \label{eq.f-f2}
     \int_{-\pi}^\pi\frac{1}{|A_{\delta_n}(e^{-i\lambda})|^2}  d\lambda = {\mathcal O}(\delta_n^{-1})  \ \ \ \mbox{and} \ \ \   \int_{-\pi}^\pi f_{\delta_n}^2(\lambda) d\lambda = {\mathcal O}(\delta_n^{-3}).
 \end{equation}
Using Lemma~\ref{le.Per}  we then get  
\begin{align*}
    \var \Big(\frac{\hat{\gamma}_{\delta_n}({0})}{\gamma_{\delta_n}({0})}\Big)
     &= \frac{1}{\gamma_{\delta_n}^2({0})} \frac{4\pi^2}{n^2}\sum_{\lambda_j \in {\mathcal G}(n)}\var\big( I_{n,\delta_n} (\lambda_{j,n})\big)\\
     & \ \ \ \ + \frac{1}{\gamma_{\delta_n}^2({0})} \frac{4\pi^2}{n^2}\sum_{{\lambda_{k,n}} \in {\mathcal G}(n)}\sum_{\substack{{\lambda_{m,n}} \in \mathcal{G}(n)\\ |m| \neq |k|}} 
     \cov \big(I_{n,\delta_n} ({\lambda_{k,n}}), I_{n,\delta_n}({\lambda_{m,n}})\big)\\
     &= \frac{1}{\gamma_{\delta_n}^2({0})} \frac{2\pi}{n}  \int_{-\pi}^{\pi} f^2_{\delta_n}(\lambda)d\lambda 
     + {\mathcal O}\Big( \frac{1}{n^2} \int_{-\pi}^{\pi} \frac{1}{ |A_{\delta_n}(e^{-i\lambda})|^2 } d\lambda \Big)  
    +  \mathcal{O}(\frac{1}{n^3\delta_n^2})\\
     & \ \ \ \ + \frac{1}{\gamma_{\delta_n}^2({0})} \Big\{ \int_{-\pi}^{\pi}\int_{-\pi}^{\pi} \frac{2\pi}{n|A_{\delta_n}(e^{-i\lambda_1})|^2 |A_{\delta_n}(e^{-i\lambda_2})|^{2}}  f_{4,U}(\lambda_{1}, -\lambda_{2}, \lambda_{2}) d\lambda_1 d\lambda_2 \\
     & \ \ \ \  + \mathcal{O}(\frac{1}{n^2\delta_n^4})\Big\}\\
     &={\mathcal O}\Big(\frac{1}{n\delta_n} + \frac{1}{n^2\delta_n} + \frac{1}{n^3\delta_n^2}\Big)  + {\mathcal O}\Big(\frac{1}{n} + \frac{1}{n^2\delta^2_n} \Big) \\
     & = {\mathcal O}(n^{-1}\delta_n^{-1}) \rightarrow 0,  
\end{align*}
by Assumption~\ref{as.3}.
\hfill $\Box$
\vspace*{0.2cm}

Before proceeding  with  the proof of Theorem~\ref{th.1} we establish the following useful  lemma the proof of which is deferred  to  the end of this section.

\vspace*{0.2cm}

\begin{lemma} \label{le.last}
The following assertions  hold true for $h\in \N$:
\begin{enumerate}
\item[(i)] $\displaystyle  \int_{-\pi}^\pi \big(\cos(\lambda h) -\rho_{\delta_n}(h)\big)^2 f^2_{\delta_n}(\lambda)d\lambda = {\mathcal O}(\delta_n^{-1}).$
\item[(ii)]$ \displaystyle
\int_{-\pi}^\pi 
\int_{-\pi}^\pi  
\big(\cos( \lambda_{1}h) -\rho_{\delta_n}(h)\big)  \big(\cos(\lambda_{2}h)-\rho_{\delta_n}(h)\big) \frac{\displaystyle 2\pi f_{4,U}(\lambda_{1}, -\lambda_{2}, \lambda_{2})}{\displaystyle n|A_{\delta_n}(e^{-i\lambda_{1}})|^2 |A_{\delta_n}(e^{-i\lambda_{2}})|^{2}} d\lambda_1d\lambda_2  =  {\mathcal O}(n^{-1})
$
\end{enumerate}
\end{lemma}
\vspace*{0.5cm}

\noindent {\bf Proof of Theorem~\ref{th.1}:}
We first  consider 
$ \hat{\rho}_{\delta_n}({h}) - \rho_{\delta_n}({h}) $. We have 
\begin{align*}
\hat{\rho}_{\delta_n}({h})-  \rho_{\delta_n}({h})\  &= \frac{\hat{\gamma}_{\delta_n}({h})}{\hat{\gamma}_{\delta_n}({0})} -\rho_{\delta_n}({h})\\
\ \ \ &=  \frac{1}{\hat{\gamma}_{\delta_n}({0})}(\hat{\gamma}_{\delta_n}({h}) - \rho_{\delta_n}({h})\hat{\gamma}_{\delta_n}({0}))\\
\ \ \ &=  \frac{1}{\frac{2 \pi}{n}\sum_{j \in \mathcal{G}(n)} I_{n,\delta_n}(\lambda_{j,n})}\Big\{\frac{2 \pi}{n}\sum_{\lambda_{j,n} \in \mathcal{G}(n)}I_{n,\delta_n}(\lambda_{j,n})\cos(\lambda_{j,n} h) \\
& \ \ \ \ - \frac{2 \pi}{n}\sum_{j \in \mathcal{G}(n)}I_{n,\delta_n}(\lambda_{j,n}) \rho_{\delta_n}({h})\Big\}\\
\ \ \ & = \frac{\gamma_{\delta_n}({0})}{\hat{\gamma}_{\delta_n}({0})} \cdot \frac{1}{\gamma_{\delta_n}({0})}{\frac{2 \pi}{n}\sum_{j\in \mathcal{G}(n)}(\cos(\lambda_{j,n} h)- \rho_{\delta_n}(h))\ I_{n,\delta_n}(\lambda_{j,n})}.
\end{align*}
From Lemma~\ref{le.gamma0}, we have that $\gamma_{\delta_n}(0)\big/\hat{\gamma}_{\delta_n}(0)  \stackrel{P}{\rightarrow} 1$.  Consider  
\begin{equation*}
   R_n:= \frac{1}{\gamma_{\delta_n}({0})}{\frac{2 \pi}{n}\sum_{j \in \mathcal{G}(n)}(\cos(\lambda_{j,n}h)- \rho_{\delta_n}(h))\ I_{n,\delta_n}(\lambda_{j,n})}.
\end{equation*}

\begin{align*}
 \var ( R_n ) 
& = \frac{1}{\gamma_{\delta_n}^2({0})}{\frac{4\pi^2}{n^2}  
\sum_{j \in \mathcal{G}(n)}\big(\cos(\lambda_{j,n}h)-\rho_{\delta_n}(h))^2\ \var(I_{n,\delta_n}(\lambda_{j,n})}\big) \\
& 
 \ \ \ \ + \frac{1}{\gamma_{\delta_n}^2(0)}\frac{4\pi^2}{n^2} 
\sum_{j\in \mathcal{G}(n)} \sum_{\substack{\lambda_{k,n}\in \mathcal{G}(n)\\ |j| \neq |k|}} 
\big(\cos(\lambda_{j,n}h)-\rho_{\delta_n}(h) \big)\big(\cos(\lambda_{k,n} h)-\rho_{\delta_n}(h)\big)\\
& \ \ \  \ \ \ \times \cov (I_{n,\delta_n}(\lambda_{j,n}), I_{n,\delta_n}(\lambda_{k,n})\big)\\
&=\frac{1}{\gamma_{\delta_n}^2({0})}\frac{4\pi^2}{n^2} 
\sum_{j \in \mathcal{G}(n)}  (\cos(\lambda_j h)- \rho_{\delta_n}(h))^2  \Big(f^2_{\delta_n}(\lambda_j)  \\
& \ \ \ \ \ \ +  \mathcal{O}\big( |A_{\delta_n}(e^{-i\lambda_{j,n}})|^{-2} n^{-1}\delta_n^{-2} + n^{-2}\delta_{n}^{-4}\Big)    \\
& \ \ \ \ + \frac{1}{\gamma_{\delta_n}^2(0)}\frac{4\pi^2}{n^2} 
\sum_{j\in \mathcal{G}(n)} \sum_{\substack{k \in \mathcal{G}(n)\\ |j| \neq |k|}}\big(\cos({\lambda_{j,n}}h)-\rho_{\delta_n}(h)\big) \big(\cos({\lambda_{k,n}}h)-\rho_{\delta_n}(h)\big)\\
& \ \ \ \ \ \  \times \Big(\frac{2\pi f_{4,U}(\lambda_{j,n}, -\lambda_{k,n}, \lambda_{k,n})}{n|A_{\delta_n}(e^{-i\lambda_{j,n}})|^2 |A_{\delta_n}(e^{-i\lambda_{k,n}})|^{2}}  + \ \mathcal{O}(n^{-2}\delta_n^{-4})\Big) \\
&=\frac{1}{\gamma_{\delta_n}^2({0})}\frac{4\pi^2}{n^2} 
\sum_{j \in \mathcal{G}(n)}  (\cos(\lambda_j h)- \rho_{\delta_n}(h))^2 f^2_{\delta_n}(\lambda_j) + {\mathcal O}\big(n^{-2}\delta_n^{-1} + n^{-3}\delta_n^{-2} \big)  \\
& \ \ \ \ + \frac{1}{\gamma_{\delta_n}^2(0)}\frac{4\pi^2}{n^2} 
\sum_{j\in \mathcal{G}(n)} \sum_{\substack{k \in \mathcal{G}(n)\\ |j| \neq |k|}}\big(\cos({\lambda_{j,n}}h)-\rho_{\delta_n}(h)\big) \big(\cos({\lambda_{k,n}}h)-\rho_{\delta_n}(h)\big)\\
& \ \ \ \ \ \  \times \frac{2\pi f_{4,U}(\lambda_{j,n}, -\lambda_{k,n}, \lambda_{k,n})}{n|A_{\delta_n}(e^{-i\lambda_{j,n}})|^2 |A_{\delta_n}(e^{-i\lambda_{k,n}})|^{2}}  + \ \mathcal{O}\big(n^{-2}\delta_n^{-2}\big) \\
&=\frac{1}{\gamma_{\delta_n}^2({0})}\frac{2\pi}{n} 
\int_{-\pi}^\pi  (\cos(\lambda h)- \rho_{\delta_n}(h))^2  f^2_{\delta_n}(\lambda)d\lambda  + {\mathcal O}\big( n^{-2}\delta_n^{-1} +   n^{-2}\delta_n^{-1} + n^{-3}\delta_n^{-2} \big) \\
& \ \ \ \ + \frac{1}{\gamma_{\delta_n}^2(0)}\int_{-\pi}^\pi 
\int_{-\pi}^\pi  
\big(\cos(\lambda_1h)-\rho_{\delta_n}(h)\big) \big(\cos(\lambda_{2}h)-\rho_{\delta_n}(h)\big)\\
& \ \ \ \ \ \  \times \frac{2\pi f_{4,U}(\lambda_{1}, -\lambda_{2}, \lambda_{2})}{n|A_{\delta_n}(e^{-i\lambda_{1}})|^2 |A_{\delta_n}(e^{-i\lambda_{2}})|^{2}} d\lambda_1d\lambda_2  + \ \mathcal{O}\big(n^{-2}\delta_n^{-1}+ n^{-2}\delta_n^{-2}\big) \\
& = R_{1,n} +R_{2,n} + {\mathcal O}(n^{-2}\delta_n^{-2} ), 
\end{align*}
with an obvious notation for $ R_{1,n}$ and $ R_{2,n}$. Note 
that  the equality before the last one holds true since the error made in the corresponding terms 
by replacing the  
Riemann sums over the grid $ {\mathcal G}(n)$ with the  integrals is of the order  $\mathcal{O}\big(n^{-2}\delta_n^{-1} )$. To elaborate,  note  the  error 
\[ \big|\frac{2\pi}{n}\sum_{j\in{\mathcal G}(n)}\big( \cos(\lambda_{j,n})-\rho_{\delta_n}(1)\big)^2 f^2_{\delta_{n}}(\lambda_{j,n}) - \int_{-\pi}^\pi \underbrace{\big(\cos(\lambda) -\rho_{\delta_n}(1)\big)^2 f^2_{\delta_n}(\lambda)}_{=:g(\lambda)}d\lambda\big|\]
 is bounded  by $2\pi \max_\lambda | g^{'}(\lambda)  |/n$. The important part in  evaluating  $ \max_\lambda | g^{'}(\lambda)  |$  is the term
\[ \frac{\big(\cos(\lambda) -\rho_{\delta_n}(1)\big)^2}{\big(\delta_n^2 +2(1-\delta_n)(1-\cos(\lambda -\lambda_0) )\big)^2} = \frac{\big(\cos(\lambda_0+t) -\rho_{\delta_n}(1)\big)^2}{ \big(\delta_n^2 +2(1-\delta_n)(1-\cos(t) )\big)^2} =:\frac{h^2_1(t)}{h^2_2(t)} =:h(t),\]
where the first equality follows using the substitution $ \lambda-\lambda_0=t$.  Differentiating $ h(t)$   and using the decomposition 
\[ \cos(\lambda_0+t) -\rho_{\delta_n}(1) = \big(\cos(\lambda_0+t) -\cos(\lambda_0)\big) + \big( \cos(\lambda_0) -\rho_{\delta_n}(1)\big),\]
we can split the derivative $ h^{'}(t)$ in the four terms 
\begin{align*}
T_1 & =\frac{\big(\cos(\lambda_0+t) - \cos(\lambda_0)\big)}{h^2_2(t)}, \ \ \ \ \    T_3 =\frac{\big(\cos(\lambda_0+t) - \cos(\lambda_0)\big)^2\sin(t)}{h^3_2(t)}, \\
T_2 & =\frac{\big(\cos(\lambda_0) - \rho_{\delta_n}(1)\big)}{h^2_2(t)}, \ \ \   \ \  T_4 =\frac{\big(\cos(\lambda_0) - \rho_{\delta_n}(1)\big)^2\sin(t)}{h^3_2(t)}. 
\end{align*}
Obtaining the maximum of the   terms $T_j$, $j=1,\ldots, 4$ via  differentiation we get  using the approximations,  $ \sin(t) \approx t$, $\cos(t) \approx 1-t^2/2$ and $ \cos(\lambda_0+t) -\cos(\lambda_0) \approx -t\sin(\lambda_0)$, after some algebra that  all four terms  are of the same order  $ {\mathcal O}(\delta_n^{-3}) $. This  implies, 
\[ \big|\frac{2\pi}{n}\sum_{j\in{\mathcal G}(n)}\big( \cos(\lambda_{j,n})-\rho_{\delta_n}(1)\big)^2 f^2_{\delta_{n}}(\lambda_{j,n}) - \int_{-\pi}^\pi \big(\cos(\lambda) -\rho_{\delta_n}(1)\big)^2 f^2_{\delta_n}(\lambda)d\lambda\big| = {\mathcal O}( n^{-1}\delta_n^{-3})\]
and therefore, 
\begin{align*}
\big| \frac{1}{\gamma_{\delta_n}^2({0})}\frac{4\pi^2}{n^2} 
& \sum_{j \in \mathcal{G}(n)}  (\cos(\lambda_j h)- \rho_{\delta_n}(h))^2 f^2_{\delta_n}(\lambda_j) - \frac{1}{\gamma_{\delta_n}^2({0})}\frac{2\pi}{n} 
\int_{-\pi}^\pi  (\cos(\lambda h)- \rho_{\delta_n}(h))^2  f^2_{\delta_n}(\lambda)d\lambda\big| \\
& 
 =  {\mathcal O}\big(n^{-2}\delta_n^{-1}\big). 
\end{align*}  

 Consider next the terms $ R_{1,n}$ and $R_{2,n}$. Using  Lemma~\ref{le.last}(i)  we get that
\[ R_{1,n}=\frac{1}{\gamma_{\delta_n}^2({0})}\frac{2\pi}{n} 
\int_{-\pi}^\pi  (\cos(\lambda h)- \rho_{\delta_n}(h))^2  f^2_{\delta_n}(\lambda)d\lambda = {\mathcal O}(n^{-1} \delta_n )\]
and
\begin{align*}
R_{2,n} & = \frac{1}{\gamma_{\delta_n}^2(0)}\int_{-\pi}^\pi 
\int_{-\pi}^\pi  
\big(\cos(\lambda_1h)-\rho_{\delta_n}(h)\big) \big(\cos(\lambda_{2}h)-\rho_{\delta_n}(h)\big)\\
& \ \ \ \ \ \  \times \frac{2\pi f_{4,U}(\lambda_{1}, -\lambda_{2}, \lambda_{2})}{n|A_{\delta_n}(e^{-i\lambda_{1}})|^2 |A_{\delta_n}(e^{-i\lambda_{2}})|^{2}} d\lambda_1d\lambda_2  \\
& = {\mathcal O}(n^{-1}\delta_n^2 ).
\end{align*}
Hence,   $ {\rm Var}(R_{n}) = {\mathcal O}\big(n^{-1}\delta_n + n^{-2}\delta_n^{-2} \big) $, that is, 
\[ \hat{\rho}_{\delta_n}({h})-  \rho_{\delta_n}({h})  = {\mathcal O}_{P}\big(n^{-(1+\alpha)/2}\big), \]
for $\alpha \in (0,1/3)$ which ensures that   $ 1+\alpha <  2-2\alpha$. 

Assertion (i) of the theorem states that  the above convergence rate of the autocorrelation estimators $ \widehat{\rho}_{\delta_n}(1)$ and $\widehat{\rho}_{\delta_n}(2)$   is transmitted to the parameter estimators 
\begin{align*}
\widehat{\alpha}_{1,\delta_n} = \frac{\hat{\rho}_{\delta_n}(1) (1-\hat{\rho}_{\delta_n}(2))}{1-[\hat{\rho}_{\delta_n}(1)]^2} & \ \ \  \mbox{and} \ \ \ 
\widehat{\alpha}_{2,\delta_n} = \frac{\hat{\rho}_{\delta_n}(2)-[\hat{\rho}_{\delta_n}(1)]^2}{1-[\hat{\rho}_{\delta_n}(1)]^2}. 
\end{align*}
This follows by a  linearization argument  and  by verifying that  the vectors of partial derivatives  of both  functions  $ g_1, g_2: (-1,1)\times(-1,1) \rightarrow \R$ with 
$ g_1(r_1,r_2) = r_1(1-r_2)/(1-r_1^2)$ and 
$ g_2(r_1,r_2)=(r_2-r_1^2)/(1-r_1^2)$,  are different from the zero  vector  $ (0,0)$ for every  $ r_2 \in(-1,1)$.
Assertion (ii) of the theorem,  follows by a   linearization argument again, the observation    that   $ \widehat{\lambda}_0=g(\widehat{a}_{1,\delta_n}, \widehat{a}_{2,\delta_n})$, where
$ g: [-2,2]\times (-1,0) \rightarrow \R$ with 
$ g(x_1,x_2) =  \arccos\big( x_1(1-x_2)/(-4x_2)\big)$  and  the fact that  for $ u =x_1(1-x_2)/(-4x_2)$, 
$$ \big( \frac{\partial g(x_1,x_2)}{\partial x_1}= (1-x_2)/\big(4x_2\sqrt{1-u^2}\big), \ \ \frac{\partial g(x_1,x_2)}{\partial x_2}= -4x_1/\big(16x_2^2 \sqrt{1-u^2}\big)\big)^\top\neq (0,0)^\top, $$
  for   every  $ x_2 \in (-1,0)$. 
\hfill $\Box$

\vspace*{0.2cm}

\noindent {\bf Proof of Lemma~\ref{le.last}:}
Consider (i).  Note that 
\[ |A_\delta(e^{-i(\lambda-\lambda_0)}|^2 = \delta^2\Big( 1 + 2\frac{1-\delta}{\delta} (1-\cos(\lambda-\lambda_0)\Big),\]
and that it suffices to consider the integral
\[ \int_0^\pi \frac{\big(\cos(\lambda h) - \rho_{\delta_n}(h)\big)^2}{ \displaystyle \delta_n^4\Big(1  + 2 \frac{1-\delta_n}{\delta_n^2}(1-\cos(\lambda-\lambda_0)\Big)^2} d\lambda. \]
We first split this integral in  two parts,  the one  for  $ |\lambda-\lambda_0| \in [0, \delta_n^r]$ and  the other for $ |\lambda-\lambda_0| \in (\delta_n^r,\pi]$ for some $ r\in (0,1]$. For the first part we get 
\begin{align*}
\int_{\lambda_0 -\delta_n^r}^{\lambda_0+\delta_n^r} & \frac{\displaystyle  \big(\cos(\lambda h) \pm \cos(\lambda_0h)-   \rho_{\delta_n}(h)\big)^2} {\displaystyle 
\delta_n^4\Big(1  + 2 \frac{1-\delta_n}{\delta_n^2}(1-\cos(\lambda-\lambda_0)\Big)^2} d\lambda\\
&  \leq  {\mathcal O}(\delta_n^{2r-4}  + \delta_n^{-2}) \int_{\lambda_0 -\delta_n^r}^{\lambda_0+\delta_n^r}\big(1 +2 \frac{1-\delta_n}{\delta_n^2}(1-\cos(\lambda-\lambda_0) \big)^2d\lambda\\
& \leq   {\mathcal O}(\delta_n^{2r-4} ) \int_{0}^{\delta_n^r}\Big(1 +\frac{1}{3}\cdot  \frac{1-\delta_n}{\delta_n^2} t^2 \Big)^{-2}dt \\
& ={\mathcal O}(\delta_n^{2r-3}).
\end{align*} 
Consider next  the second part of the integral. We have,
\begin{align*}
\int_{|\lambda-\lambda_0 |> \delta_n^r} & \frac{\displaystyle  \big(\cos(\lambda h) \pm \cos(\lambda_0h)-   \rho_{\delta_n}(h)\big)^2} {\displaystyle 
\delta_n^4\Big(1  + 2 \frac{1-\delta_n}{\delta_n^2}(1-\cos(\lambda-\lambda_0)\Big)^2} d\lambda\\
&\leq {\mathcal O}(\delta_n^{-4}) \int_{|\lambda-\lambda_0 |> \delta_n^r} \frac{(\lambda-\lambda_0)^2 +\delta_n^2}{\displaystyle \Big(1 + 2 \frac{1-\delta_n}{\delta_n^2} (1-\cos(\lambda-\lambda_0)\Big)^2}d\lambda.
\end{align*} 
After some algebra we get that 
\[  \int_{|\lambda-\lambda_0 |> \delta_n^r} \frac{(\lambda-\lambda_0)^2}{\displaystyle \Big(1 + 2 \frac{1-\delta_n}{\delta_n^2} (1-\cos(\lambda-\lambda_0)\Big)^2}d\lambda= {\mathcal  O}(\delta_n^3),\]
and 
\[ \int_{|\lambda-\lambda_0 |> \delta_n^r} \frac{ \delta_n^2}{\displaystyle \Big(1 + 2 \frac{1-\delta_n}{\delta_n^2} (1-\cos(\lambda-\lambda_0)\Big)^2}d\lambda = {\mathcal O}(\delta_n^3).\]
That is, 
\[ \int_{|\lambda-\lambda_0 |> \delta_n^r}  \frac{\displaystyle  \big(\cos(\lambda h) \pm \cos(\lambda_0h)-   \rho_{\delta_n}(h)\big)^2} {\displaystyle 
\delta_n^4\Big(1  + 2 \frac{1-\delta_n}{\delta_n^2}(1-\cos(\lambda-\lambda_0)\Big)^2} d\lambda =  {\mathcal O}(\delta_n^{-1}).\]
Hence,
\[ \int_{0}^\pi \big(\cos(\lambda h) - \rho_{\delta_n}(h)\big)^2 f_{\delta_n}^2(\lambda) d\lambda = {\mathcal O}\big(\delta_n^{2r-3} + \delta_n^{-1}\big) = {\mathcal O}(\delta_n^{-1}),\]
for $ r=1$. 

Consider (ii).  Since $ f_{4,U}$ is bounded  it suffices to evaluate the integral
\[ T_n:= \frac{1}{n}\Big(\int_{-\pi}^\pi \frac{(\cos(\lambda h) -\rho_{\delta_n}(h))}{|A_{\delta_n}(e^{-i (\lambda-\lambda_0)})|^2 }d\lambda \Big)^2.\]
Arguing as in the proof of (i)  and using
\[ \int_0^\pi \frac{1}{1+at^2}dt = {\rm arctan} (\sqrt{a} t)/\sqrt{a} \Big|^\pi_0 \ \ \mbox{and} \int_0^\pi \frac{t }{1+at^2}dt = (2a)^{-1} \log(1+at^2)  \Big|^\pi_0,\]    
we get  after some algebra that  $ T_n= {\mathcal O} (n^{-1} \delta_n^{-1} \log(n) + n^{-1}) ={\mathcal O}(n^{-1})$. 
\hfill $\Box$

 \end{document}